\documentclass[11pt, twoside, a4paper]{article}
\usepackage{graphicx, graphics,latexsym,amssymb, amsmath}
\usepackage[mathscr]{euscript}
\usepackage{url}
\usepackage[english]{babel}

\title{
Smooth Riemannian Structures on Dessins d'Enfants}
\author{Jean-Marie Morvan}
\newcommand{\Addresses}{{
  \bigskip
  \footnotesize

  Jean-Marie Morvan,
\textsc{Universit\'{e} de Lyon, CNRS UMR $5208$,
Universit\'{e} Claude Bernard Lyon $1$, Institut Camille Jordan,
 $43$ blvd du $11$ Novembre $1918$, $F$-$69622$ Villeurbanne-Cedex, France}
 {\rm and}
 \textsc{King Abdullah University of Science and Technology, $C.E.M.S.E$,  Thuwal $23955$-$6900$, Saudi Arabia}\\
\\\par\nopagebreak
  \textit{E-mail address}~: \texttt{morvan@math.univ-lyon1.fr} {\rm and}
\texttt{Jean-Marie.Morvan@KAUST.EDU.SA}
}}

\def\cal{\mathcal}
\font\bb=msbm10
\def\Z{\hbox{\bb Z}}
\def\R{\hbox{\bb R}}
\def\E{\hbox{\bb E}}
\def\S{\hbox{\bb S}}
\def\N{\hbox{\bb N}}

\def\C{\hbox{\bb C}}

\def\H{\hbox{\bb H}}
\def\D{\hbox{\bb D}}

\newtheorem{lemma}{Lemma}
\newtheorem{theorem}{Theorem}
\newtheorem{proposition}{Proposition}
\newtheorem{definition}{Definition}
\newtheorem{corollary}{Corollary}

\newenvironment{Myitemize}{\begin{itemize}
}{\end{itemize}}

\begin{document}

\maketitle

\begin{abstract}{\footnotesize
We show how to define  a canonical Riemannian metric on a {\it "dessin d'enfants"} drawn on a topological surface.
This gives a possible explanation of a claim of A. Grothendieck \cite{grothendieck1984esquisse}.}
\end{abstract}

\vspace{20pt}

In his famous paper {\it Esquisse d'un programme} \cite{grothendieck1984esquisse},  Alexandre  Grothendieck  introduced a new view point in the study of maps on surfaces (that he called {\it dessins d'enfants}). At the end of this text, insisting on the interest to pursue his ideas, he  briefly pointed out a link with Riemannian geometry:

\vspace{10pt}

\noindent {\it "Depuis $1977$, dans toutes les questions (comme dans ces deux derniers
th\`{e}mes que je viens d'\'{e}voquer) o\`{u} interviennent des cartes bidimensionelles,
la possibilit\'{e} de les r\'{e}aliser canoniquement sur une surface conforme, donc
sur une courbe alg\'{e}brique complexe dans le cas orient\'{e} compact, reste
en filigrane constant dans ma r\'{e}flexion.  Dans pratiquement tous les cas
(en fait, tous les cas sauf celui de certaines cartes sph\'{e}riques avec "peu
d'automor-phismes") une telle r\'{e}alisation conforme implique en fait une
m\'{e}trique riemanienne canonique,
ou du moins, canonique \`{a} une constante
multiplicative pr\`{e}s."
}
\vspace{10pt}

\noindent {\it ["Since $1977$, in all the questions (such as the two last themes evoked above) where two-dimensional maps occur, the possibility of realising them canonically on a conformal surface, so on a complex algebraic curve in the compact oriented case, remains constantly in filigree throughout my reflection. In practically every case (in fact, in all cases except that on certain spherical maps with "few automorphisms") such a conformal realisation implies in fact a canonical Riemannian metric, or at least, canonical up to a multiplicative constant."]}

\vspace{10pt}

\noindent Two statements in these lines could be clarified~: First of all, A.  Grothendieck does not specify the regularity of the canonical metric. Moreover,  the notion of   {\it "map with few automorphisms"} is not precisely defined.
\begin{Myitemize}
\item
{\it A priori}, to find a Riemannian metric canonically associated to a dessin ${\cal D}$ on a surface $X$, one could build  a triangulation associated to the dessin,  endow this triangulation with a piecewise linear structure, and then, build the (singular) Euclidean metric defined by affecting the length $1$ to each edge. Such a metric is {\it not smooth} in general (only continuous, since it  may have conical singularities at the vertices).  This construction does not need any restriction on the automorphism group of the dessin.

\item
These considerations lead to look for a {\it smooth} canonical metric associated to  the dessin ${\cal D}$. The natural way is to  build the standard conformal (that is, complex) structure associated to ${\cal D}$, and then to consider  Riemannian metrics invariant by the group ${\rm Aut}(X)$ of biholomorphisms of $X$. The Poincare-Klein-Koebe uniformisation theorem gives a Riemannian metric with constant Gaussian curvature   in  the conformal class of Riemannian metrics defined on $X$. This Riemannian metric is unique if the genus of $X$ is strictly greater than $1$, unique up to a scaling constant if the genus is $1$, and invariant by the group ${\rm Aut}(X)$ of biholomorphisms of $X$. However, this construction fails if the genus of $X$ is $0$. In this situation, $X$ is biholomorphic to the Riemann sphere $\hat{\C}$, the modular space of $X$ is reduced to a point and a dessin does not induce any information on the conformal structure of $X$. Moreover, the standard metric with constant Gaussian curvature $1$ in the Riemann sphere is not invariant by the group ${\rm Aut}(\hat{\C})$ of biholomorphisms of $\hat{\C}$.  In order to find a canonical metric associated to a dessin drawn on $\hat{\C}$, the idea is to replace the (too large) group ${\rm Aut}(\hat{\C})$ by the (finite) subgroup of biholomorphisms ${\rm Aut}(\hat{\C},{\cal D})$ that preserve ${\cal D}$ and to build a Riemannian metric invariant only by ${\rm Aut}(\hat{\C},{\cal D})$.  We propose two different constructions.
 \begin{itemize}
 \item
 The first one is based on the average (over  ${\rm Aut}(\hat{\C},{\cal D})$) of the metrics obtained by pullback of the standard metric of the round $2$-sphere of radius $1$. We can build this metric without any restrictions on  ${\rm Aut}(\hat{\C},{\cal D})$. However, its Gaussian curvature is not constant in general.
 \item
 The second one mimics the hyperbolic situation, but excludes the case where ${\rm Aut}(\hat{\C},{\cal D})$ is cyclic. The metric we get is invariant by  ${\rm Aut}(\hat{\C},{\cal D})$ and of constant Gaussian curvature $1$.     The opinion of the author is that it  may be the one considered by A. Grothendieck, who excluded in his text, {\it "certain spherical maps with "few automorphisms""}, although he did not mention that he looked for a metric with constant Gaussian curvature.
 \end{itemize}
 In both case, these  metrics coincide with the standard metric of the round sphere $\S^2$ of radius $1$ when ${\rm Aut}(\hat{\C},{\cal D})$ is  a subgroup of $SO(3)$.
\end{Myitemize}

\vspace{10pt}

\noindent Let us now state the theorem corresponding to the second construction (the (non smooth) continuous situation is mentioned in section \ref{CANRI}, and the result corresponding to the first construction is stated in section \ref{SFSIT}), see notations and definitions below~:

\begin{theorem} \label{THPR}
Let $(X, {\cal D})$ be a dessin on a (closed oriented) topological surface.
\begin{Myitemize}
\item
 If ${\bf g}_X > 1$, then $(X, {\cal D})$ admits a canonical Riemannian metric with constant Gaussian curvature $-1$ induced by ${\cal D}$, invariant by ${\rm Aut}(X)$.
\item
If ${\bf g}_X =1$,  then $(X, {\cal D})$ admits a canonical  Riemannian  flat  metric induced by ${\cal D}$,  unique up to a scaling constant, invariant by ${\rm Aut}(X)$.
\item
If ${\bf g}_X =0$ and  ${\rm Aut}(X,{\cal D})$ is not cyclic, then $(X, {\cal D})$ admits a canonical  Riemannian metric of constant Gaussian curvature $1$ invariant by  ${\rm Aut}(X,{\cal D})$. In particular,  if ${\rm Aut}(X,{\cal D})$ is  a subgroup of $SO(3)$ (canonically embedded in $PSL(2,\C)$), this metric coincides with the standard metric of the round $2$-sphere of radius $1$.
 \end{Myitemize}
\end{theorem}

\noindent  This text trying to be as self contained as possible, we present in the following paragraphs some backgrounds on conformal geometry, Riemannian geometry and the theory of dessins d'enfants. Many results of this paper are simple reminders (in particular in Riemannian geometry, conformal geometry and complex analysis) fixing the notations and giving the essential results needed for the proof of Theorem \ref{THPR}. We refer to classical or more recent  books as \cite{donaldson2011riemann}, \cite{REYSSAT} \cite{bost1992introduction}, \cite{guillarmoumoroianu}, \cite{earp1997}, \cite{gu2008computational}, \cite{dai2008variational} for complete and detailed studies of these topics.

\vspace{10pt}

 \noindent The author would like to thank  J. Germoni, E. Toubiana, J. Wolfart, A. Zvonkin, for useful mails, discussions, and improvements.

\section{Reminders on some classical groups}
 \subsection{The groups $GL(2,\C)$, $PGL(2,\C)$, $SL(2,\C)$, $PSL(2,\C)$}
\begin{itemize}
\item
The {\it linear group}  $GL(2,\C)$  is the group of linear automorphisms of $\C^2$, identified to the group of invertible complex $(2,2)$-matrices
$\begin{pmatrix}
a & b \\
c & d
\end{pmatrix},
$ where
 $a,b,c,d \in \C, ad-bc \neq 0$. Its center $C$ is the subgroup of homothecies identified to the subgroup of matrices
  $\begin{pmatrix}
\lambda & 0 \\0 & \lambda
\end{pmatrix},
$
where $\lambda \in \C^{*}$.

\item
The {\it projective linear group} is the quotient $PGL(2,\C)= GL(2,\C)/C$.  It can be identified to the group of invertible complex $(2,2)$-matrices
$\begin{pmatrix}
a & b \\c & d
\end{pmatrix},
$ where
 $a,b,c,d \in \C, ad-bc =1$.

\item
 The {\it special linear group}  $SL(2,\C)$  is the normal subgroup of $GL(2,\C)$ defined as the kernel of the determinant
 homomorphism~:
 $${\mbox{\rm det }}: GL(2,\C) \to  \C^*.$$
 Its center $SC$ is the subgroup of homothecies
  $\begin{pmatrix}
\lambda & 0 \\0 & \lambda
\end{pmatrix},
$
where $\lambda \in \C^{*},  |\lambda|=1$.

\item
Finally, the {\it special projective linear group}  $PSL(2,\C)$ is the quotient  $SL(2,\C)/SC$.
\end{itemize}

\noindent The following result is clear~:

\begin{proposition} \label{ISOM1}
The groups $PSL(2,\C)$ and $PGL(2,\C)$ are isomorphic. Each element of them can be represented by a (2,2)-matrix
 $\begin{pmatrix}
a & b \\c & d
\end{pmatrix},
$ where
 $a,b,c,d \in \C$, with  $ad-bc =1$.
\end{proposition}

\subsection{The group $SO(3)$ and its finite subgroups}
 Let $SO(3)$ denotes the group of positive isometries of $\R^3$. Every  element $r$ of $SO(3)$ different to the identity is a rotation of $\R^3$ with axis $\delta_r$. Such a rotation $r$ acts on $\S^2$, with two fixed points $x_r$ and $-x_r$ (the intersections of  $\delta_r$ with $\S^2$).

 \vspace{10pt}

 \noindent Let us describe the finite subgroups of $SO(3)$.
\begin{theorem}\label{SO38}
Any finite subgroup $K$ of $SO(3)$ is isomorphic to one of the following groups~: A cyclic group $\Z_n$, a diedral group ${\cal D}_n$,  $(n \in \N^{*})$, the symmetric group ${\cal S}_4$,   the alternate group ${\cal A}_4$,   the alternate group ${\cal A}_5$.
\end{theorem}

\noindent {\footnotesize {\bf Reminder~- Sketch of proof of Theorem \ref{SO38}~-}  The subgroup of rotations $r \in K$ fixing a couple $(x,-x) \in \S^2 \times \S^2$ can be identified to a subgroup of $SO(2)$ (the subgroup of rotations $r \in K$ in the plane orthogonal to the axis $\Delta_r$ of $r$). Then, it is a cyclic group $K_{x}$ isomorphic to $\Z_n$ for some $n \in \N^*$.
 Let ${\cal F}$ be the set of fixed points of any element of $K$, that is,
 $${\cal F} = \{x \in \S^2; \exists h\in K \backslash \{{\rm Id}\}, h(x)=x\}.$$
Let $x \in {\cal F}$, ({\it resp.}  $y \in {\cal F}$).  It is clear that $K_x$ and $K_y$ are conjugate~: There exists $g \in K$ such that $K_y=g^{-1}K_xg$. In particular, $|K_y|=|K_x|$. Then all cyclic groups $K_x$, $x \in K$, have the same cardinality. Let
 $${\cal P}=O_1 \cup ...\cup O_k$$
 be the partition of ${\cal F}$ into the  orbits of $K$. Since the stabilizer subgroups of any element of an orbit $O_i$ are conjugate, we can define $\epsilon_i$ as its order. Classically, the class formula and Burnside formula imply that
 $${\rm card}(O_i) = \frac{{\rm card}(K)}{\epsilon_i},$$
 and,
 $$\sum_{i=1}^k \frac{1}{i} = k-2 + \frac{2}{{\rm card} (K)}.$$
 This equation implies $k=2$ or $3$.

 \begin{Myitemize}
 \item
 If $k=2$, the only possible triplet  $(\epsilon_1, \epsilon_2; {\rm card}(K))$  is
 $$({\rm card}(K), {\rm card}(K); {\rm card}(K)).$$
 In this case, $K$ is isomorphic to $\Z_n$.

 \vspace{10pt}

 \noindent In all other cases, we remark that $K$  has three orbits.

 \item
 If $k=3$, the only possible triplets  $(\epsilon_1, \epsilon_2, \epsilon_3; {\rm card}(K))$ are the following~:
   \begin{itemize}
   \item
   $(2,2,n;2n)$,
   \item
   $(2,3,3;12)$,
   \item
   $(2,3,4;24)$,
   \item
   $(2,3,5;60)$.
   \end{itemize}
 \end{Myitemize}

\begin{enumerate}
 \item
  Let us study the case  $(2,2,n;2n)$. The group $K$ is isomorphic to the diedral group $\D_n$. The orbit $O_3$ is a subset of two fixed points. The orbits $O_1$ and $O_2$ are subsets of $n$ fixed points.

  \item
  Let us study the case $(2,3,3;12)$. 
   The group $K$ is isomorphic to the tetrahedral group ${\cal A}_4$ (it is the group of symmetry of the regular tetrahedron).  The cardinal of the orbit $O_1$ is $6$, The cardinals of the orbits $O_2$ and $O_3$ are $4$.

   \item
   An analogous study of the case $(2,3,4;24)$ shows that
   the group $K$ is isomorphic to the octahedron group ${\cal A}_4$ (it is the group of symmetry of the regular octahedron).

  \item
    An analogous study of the case  $(2,3,5;60)$ shows that the group $K$ is isomorphic to the icosahedron group ${\cal A}_5$ (it is the group of symmetry of the regular icosahedron).

  \end{enumerate}
  }

\begin{figure}[h]
\begin{minipage}[b]{0.3\linewidth}
\centering
\includegraphics[width=\textwidth]{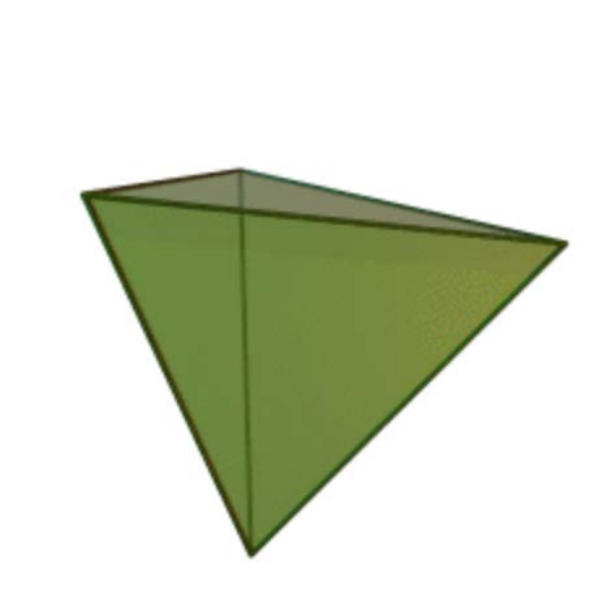}
\end{minipage}
\hspace{0.3cm}
\begin{minipage}[b]{0.3\linewidth}
\centering
\includegraphics[width=\textwidth]{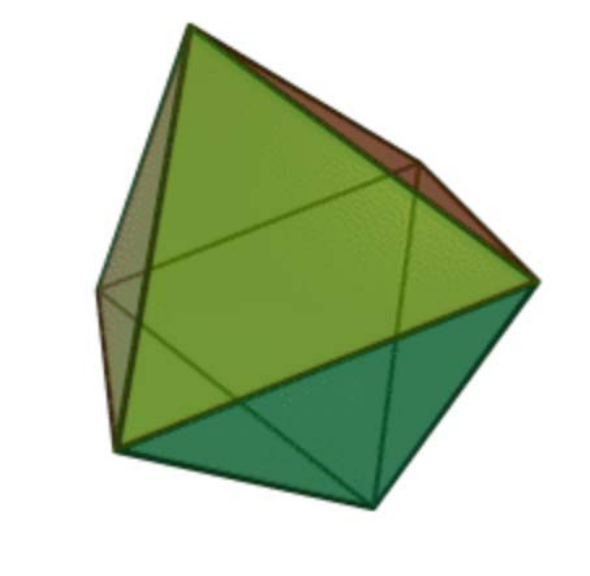}
\end{minipage}
\hspace{0.5cm}
\begin{minipage}[b]{0.3\linewidth}
\centering
\includegraphics[width=\textwidth]{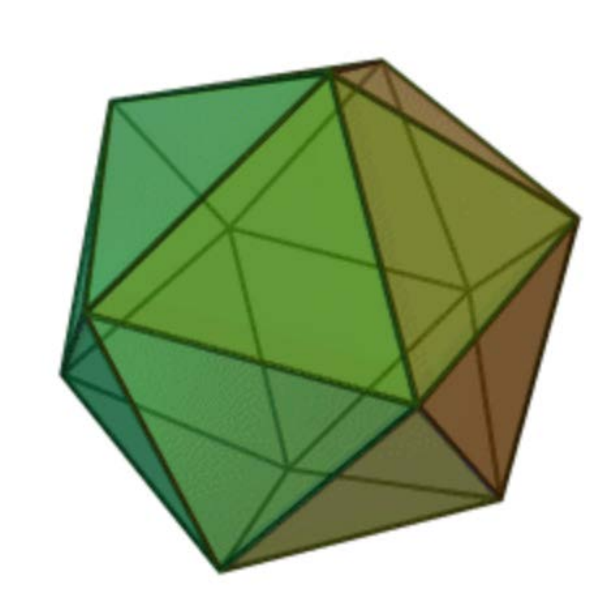}
\end{minipage}
\caption{A tetrahedron, an octahedron and an icosahedron, (From Wikipedia, Platonic solid)
 }
\label{teocic}
\end{figure}

  \noindent We remark that the groups ${\cal A}_4$, ${\cal S}_4$, ${\cal A}_5$ are the groups of symmetry of  Platonic solids (see  Figure \ref{teocic} from Wikipedia, Platonic solids).

\begin{proposition} \label{SO3}
There exists a canonical embedding of $SO(3)$ into $PSL(2,\C)$.
\end{proposition}

\noindent A classical proof of Proposition \ref{SO3} consists of using the algebra of quaternions.
A more geometrical proof can be done as follows~: {\it Via} the stereographic projection $s$ described in section \ref{RUTH6}, any rotation of the sphere $\S^2$ can be transported to a bijection from $\hat{\C}$ to itself. A direct computation shows that this bijection belongs to $PSL(2,\C)$.

\vspace{10pt}

\noindent Proposition \ref{SO3} implies that we can consider $SO(3)$ as a subgroup of  $PSL(2,\C) \simeq PGL(2,\C)$.

\subsection{The finite subgroups of $PSL(2,\C)$} \label{FINSUB}
Theorem \ref{SO38} can be extent to the finite subgroups of $PSL(2,\C)$. We identify $SO(3)$ with its image by the embedding  given in Proposition \ref{SO3}.
Let $\mu_{n}^{*}$ denote the set of primitive roots of unity in $\C$, and define the equivalence relation $\sim$ on $\mu_{n}^{*}$ by
$$z \sim z' \iff z' \in \{z, z^{-1}\}.$$

\begin{proposition}
Let $h \in PSL(2,\C)), h \neq {\rm Id}$.
\begin{enumerate}
\item
 If the order of $h$ is finite, then $h$ has two fixed points in $\hat{\C}$.

\item
    \begin{enumerate}
       \item
          The following assertions are equivalent~:
 \begin{enumerate}
 \item
 The order of $h$ is $n >1$;
 \item
 $h$ is conjugate to a homothecy  $z \to \zeta z$, where $[\zeta] \in \mu_{n}^{*}/\sim$.
 \end{enumerate}
 \item
The space  $\mu_{n}^{*}/\sim$ classify the conjugacy classes of order $n$.
\end{enumerate}
\end{enumerate}
\end{proposition}

\begin{theorem} \label{SGSO}
\begin{enumerate}
\item \label{SGSO1}
 Each finite subgroup of $PSL(2,\C)$ is isomorphic to one of the following  groups~:  A cyclic group $\Z_n$, a diedral group ${\cal D}_n$, the symmetric group ${\cal S}_4$,   the alternate group ${\cal A}_4$,   the alternate group ${\cal A}_5$.
\item \label{SGSO2}
All subgroups of one of these categories are conjugate in   $PSL(2,\C)$. In particular, they are conjugate to a subgroup of $SO(3,\R)$.
\end{enumerate}
\end{theorem}

\noindent The canonical embedding of $SO(3)$ into  $PSL(2,\C)$ and item \ref{SGSO2} of \ref{SGSO} allows to choose a particular "standard" element of $SO(3)$ to represent each  conjugacy class of a finite subgroup of $PSL(2,\C)$. This is the goal of Lemma \ref{REP5}. If $h,k \in PSL(2,\C)$, we denote by $\langle h, k\rangle$ the subgroup spent by $h$ and $k$.

\begin{lemma} \label{REP5}
\item
The following subgroups of $PSL(2,\C)$ are included in $SO(3)$~:
\begin{Myitemize}
\item
${\bf C}_n=\langle z \to \zeta z \rangle$, where $[\zeta] \in \mu_{n}^{*}/\sim$;
\item
${\bf D}_n=\langle z \to \zeta z, z \to \frac{1}{z} \rangle$, where $\zeta \in \mu_{n}^{*}/\sim$;
\item
${\bf A}_4=\langle z \to j z, z \to \frac{z+\sqrt{2}}{\sqrt{2}z-1} \rangle$;
\item
${\bf S}_4=\langle z \to \iota z, z \to \frac{z+1}{z-1} \rangle$;
\item
${\bf A}_5=\langle z \to \iota z, z \to \frac{z+ \Delta}{\Delta z-1} \rangle$, where $\delta$ is a fifth primitive root of unity, and $\Delta = \sqrt{1 -\delta - \frac{1}{\delta}}$.
\end{Myitemize}
\end{lemma}

\begin{theorem} \label{SUBB}
\begin{enumerate}
\item \label{SUBB1}
The cyclic subgroups of order $n$ of $PSL(2,\C)$ are conjugate to the subgroup ${\bf C}_n$;

\item \label{SUBB2}
The diedral groups of order $n$ are conjugate to the subgroup ${\bf D}_n$;

\item \label{SUBB3}
The subgroups  of $PSL(2,\C)$ isomorphic to ${\cal A}_4$ are conjugate to the subgroup ${\bf A}_4$;

\item \label{SUBB4}
The subgroups  of $PSL(2,\C)$ isomorphic to ${\cal S}_4$ are conjugate to the subgroup ${\bf S}_4$;

\item \label{SUBB5}
The subgroups  of $PSL(2,\C)$ isomorphic to ${\cal A}_5$ are conjugate to the subgroup ${\bf A}_5$.
\end{enumerate}
\end{theorem}

\noindent  We remind that
\begin{Myitemize}
\item
${\cal A}_4 \simeq \langle (123),(12)(34)\rangle$;

\item
${\cal S}_4 \simeq \langle (1234),(12)\rangle$;

\item
${\cal A}_5 \simeq \langle (12345),(12)(34)\rangle$.
\end{Myitemize}

\vspace{10pt}

\noindent For further use, we need to compute the normalizers of the finite subgroups of $SO(3)$. If $K$ is a finite subgroup of $SO(3)$, we denote by $N(K)$ its normalizer in $PSL(2,\C)$. We have the following result (\cite{cheltsov2016} for instance)~:

\begin{theorem} \label{NORMALIZER}
\begin{enumerate}
\item
The normalizer in  $PSL(2,\C)$ of any non cyclic finite subgroup of $SO(3)$ is finite and included in $SO(3)$ .
\item
More precisely,
\begin{itemize}
\item
for any  $n \in \N^*$, $N({\bf D}_n)={\bf D}_n$,
\item
$N({\bf A}_4)={\bf S}_4$,
\item
$N({\bf S}_4)={\bf S}_4$,
\item
$N({\bf A}_5)={\bf A}_5$.
\end{itemize}
\end{enumerate}
\end{theorem}

\noindent We remark however that the normalizer of a cyclic subgroup ${\bf C}_n$ is  isomorphic to the group of $2 \times 2$ diagonal matrices with complex coefficients. In particular, it is not finite,  nor in $SO(3)$.

\section{Conformal and holomorphic maps in the Euclidean plane}
We denote by $\E^2$ the (oriented) Euclidean space, that is, the oriented real two dimensional vector space $\R^2$ endowed with its standard scalar product $g_{\E^2}$. In the following, we identify  $\E^2$ with $\C$, endowed with the standard scalar product $g_{\C}$.

\subsection{Conformal maps}
Let $U$ and $V$ be open subsets of $\E^2$.
\begin{definition}
let $f: U \to V$ be a map.
\begin{Myitemize}
\item
The map $f$ is called {\rm conformal} if it preserves the angles, that is, if it satisfies the following property~: For all $p \in U$, any $X_p \in T_p\E^2$ and $Y_p \in T_p\E^2$,
\begin{equation} \label{ANGLE}
\angle (f_p(X_p),f_p(Y_p))=\angle (X_p,Y_p).
\end{equation}
\item
The map $f$ is called {\rm anti-conformal} if it preserves the absolute values of the angles (computed in $]-\pi, +\pi[$), and reverses the orientation.
\item
A bijective conformal map is called a {\rm conformal transformation}.
\end{Myitemize}
\end{definition}

\noindent In other words, if $g_0$ denotes the scalar product of $\E^2$, a map $f: U \to V$ is {\it conformal} if there exists a real function $\lambda$ defined on $U$ such that or all $p \in U$, for all $X_p \in T_p\E^2$ and $Y_p \in T_p\E^2$,
$$g_0(f_p(X_p), f_p(Y_p))= e^{\lambda}g_0(X_p,Y_p).$$

\noindent We can remark that a conformal transformation preserves the orientation. A map $f$ is {\it anticonformal} if and only if $\bar{f}$ is conformal. The link between conformal maps and holomorphic map is the following (we identify  $\E^2$ with  $\C$)~:

\begin{proposition}
A map $f: U \to V$ is conformal if and only if it is holomorphic and satisfies $f'(p) \neq 0$ for every $p \in U$.
\end{proposition}

\noindent In particular, $f$ is  a conformal transformation if and only if $f$ is {\it biholomorphic} (that is, $f$ and $f^{-1}$ are holomorphic).

\subsection{Some reminders on holomorphic maps} \label{REMIND}
 Let us recall the well known properties on holomorphic maps defined on a domain $D$.
 We will use the following theorems~:

 \begin{theorem} \label{POINT}
 Let $D$ be a domain of $\C$, $p \in D$.
  If $f$ is a holomorphic  function defined on  $D \backslash \{p\} \subset \C$, and bounded on a neighborhood of $p$, then $f$ can be extended as a holomorphic function on $D$.
  \end{theorem}

  \begin{theorem} \label{STRAIGHTL}
  If $f$ is a continuous map defined on a domain $D \subset \C$, holomorphic on $D$ except at most at the points of a straight line, then $f$ is holomorphic at every point of $D$.
  \end{theorem}

  \noindent In particular, if $f$ is a continuous function defined on a domain $U \subset \C$, holomorphic on $U$ except at most at a finite subset of points, then $f$ is holomorphic at every point of $U$.

  \begin{theorem}~-~{\bf The Riemann mapping Theorem~-}\label{RIEMANNTH}
  If $U$ is a non-empty simply connected open subset of $\C$, different to $\C$, then there exists a biholomorphic bijection $f$ from $U$ onto the open unit disk
  $$D=\{z\in \C, |z| <1\}.$$
  \end{theorem}

  \noindent As a consequence of Theorem \ref{RIEMANNTH}, a non-empty simply connected open subset $U$ of $\C$, different to $\C$ is also biholomorphic to the upper plane $\C^+= \{(x,y), y \geq 0\}$ of $\C$, since this upper plane is itself a simply connected open subset $U$ of $\C$ different to $\C$.
   Although Theorem \ref{RIEMANNTH} is a purely a existence theorem, the following Schwarz-Christoffel Theorem gives an explicit expression of the biholomorphism from  $U$ to the upper plane of $\C$, when $U$ is a polygonal region.

  \begin{theorem}~-~{\bf The Schwarz-Christoffel Theorem}~- \label{RIEMANNCHRIS}
  Let $P$ be  a polygonal region in $\C$, with $n$ vertices and interior angles $\alpha_1, \alpha_2,...,\alpha_n$.
  The primitive of the function
  \begin{equation} \label{SCH}
  f(z)=A(z-x_1)^{\frac{\alpha_1}{\pi}-1}(z-x_1)^{\frac{\alpha_2}{\pi}-2}... (z-x_n)^{\frac{\alpha_1}{\pi}-n}
  \end{equation}
  (where $A$ is a nonzero constant),
  maps the upper half plane $\C^+$ to $P$, in such a way that the real axis is sent on the edges of $P$, and the points $x_1, \ldots ,x_n$  on the real axis are sent on  the vertices of $P$.
  \end{theorem}

  \noindent  Of course, if the polygonal region is bounded, only $n-1$ angles are included in the formula \ref{SCH}. As an example, if $t$ is a right triangle with angle $\frac{\pi}{6}, \frac{\pi}{3}, \frac{\pi}{2}$, then a possible mapping $f$ mapping $t$ onto $\C^+$ is a primitive of the function
  $$f(z)= \frac{1}{z^{\frac{1}{2}}(z+1)^{\frac{2}{3}}}.$$

\section{Riemann surfaces}
\subsection{Definition of Riemann surfaces}
A Riemann surface is   a  {\it complex $1$-dimensional analytic manifold}. Let us be more precise~:

\begin{definition}
Let $X$ be an (connected oriented) topological surface endowed with an atlas $\{(U_{\alpha}, z_{\alpha}, V_{\alpha})\}_{\alpha  \in I}$, where $\{U_{\alpha}\}_{\alpha \in I}$ covers $X$ by open subsets and for every $\alpha$,
$$z_{\alpha} : U_{\alpha} \to V_{\alpha} \subset \E^2 \simeq \C$$ is a homeomorphism from $U_{\alpha}$ onto an open set $z_{\alpha}(U_{\alpha})$, such that, if $U_{\alpha} \cap U_{\beta} \neq \emptyset$,
$$z_{\beta} \circ z_{\alpha}^{-1} : z_{\alpha}(U_{\alpha} \cap U_{\beta}) \to z_{\beta}(U_{\alpha} \cap U_{\beta})$$
is a conformal transformation. The surface $X$ is called a {\rm Riemann surface}. One says that $X$ is endowed with a {\rm conformal structure}.
\end{definition}

\noindent In other words, a {\it Riemann surface} is a (connected) topological surface whose transition functions are conformal bijections between open subsets of $\C$. Such an atlas is called a {\it complex or conformal atlas}.  If the union of two conformal atlases on $X$ is still a conformal atlas, they are called {\it equivalent}. An equivalence class of conformal atlas is called a {\it conformal structure} on $X$. We will denote by the generic letter ${\cal C}$ the conformal structure defining a Riemann surface.

\vspace{10pt}

\noindent  Two Riemann surfaces are said to be {\it isomorphic} if there exists a biholomorphic bijection from the first one to the second one. An {\it automorphism} of a Riemann surface $X$ is an isomorphism from $X$ to itself.

\begin{definition}
Let $X_1$ and $X_2$ be Riemann surfaces. A map $f:X_1 \to X_2$  is called {\rm conformal} ({\rm resp.} {\rm holomorphic}) at a point $p \in X_1$ if there exists a chart $(U_1, \phi_1, V_1)$ around $p$, a chart $(U_2, \phi_2, V_2)$ around $f(p)$ such that $\phi^{-1}_2 \circ f \circ \phi^{1}_1$ is conformal ({\rm resp.} holomorphic).
\end{definition}

\subsection{The Riemann uniformization theorem} \label{RUTH6}
The following standard surfaces are endowed with a {\it canonical} structure of Riemann surface~:
\begin{Myitemize}
\item
$\C$ (it is obvious~!)

\item
Any connected open set of $\C$ and in particular  the {\it hyperbolic plane}
 $$\H=\{z \in \C; \mbox{im}(z) >0\},$$
 that is, the upper half plane of $\C$, isomorphic (as a Riemann surface) to the unit (open) disc $\D=\{z \in \C; |z|<1\}$;

\item
The {\it Riemann sphere}  $\hat{\C}$ whose underlying set is $\C \cup \infty$ endowed with an atlas of two charts~: $(U_1=\C,\varphi_1, V_1=\C), (U_2=\C^{*} \cup \{\infty\}, \varphi_2, V_2=\C^{*} \cup \{\infty\})$, where
$\varphi_1 = id$ and $\varphi_2$ is defined by
$$\begin{cases}
\varphi_2 (z) = \frac{1}{z}, \mbox{ {\rm if} } z \in \C\\
\varphi_2 (\infty)= \infty
\end{cases}
$$
The transition function is the function $$\psi : \C^* \to \C^*,$$
 defined by $\psi(z)=\frac{1}{z}$.
The Riemann surface $\hat{\C}$ is a topological $2$-sphere (in particular it is connected and compact), as it is easily shown by the {\it stereographic projection}
 that we describe now~:
One identifies $\R^3$ with $\C \times \R$, we denote by $\S^2$ its unit sphere, $\C$ being the "horizontal plane". The north pole $n$ is the point $(0,1)$. The {\rm stereographic projection} is the map
$$s: \S^2 \subset \C \times \R \to \hat{\C},$$
defined for every point $m$ on $\S^2$ by
$$
\begin{cases}
s(m) &= nm \cap \C, \mbox{ {\rm if} } m\neq n,\\
s(n) &= \infty,
\end{cases}
$$
where $nm$ denotes the line throwing $n$ and $m$. The map $s$ is an homeomorphism defined as follows~:
 for all $(z,t) \in \C \times \R$ such that $|z|^2 + t^2=1$,
 $$s(z,t) = \frac{z}{1-t}.$$
 Its inverse
 $$s^{-1} : \hat{\C} \to \S^2 \subset \C \times \R$$ satisfies~:
 for all $z\in \C$,
$$s^{-1}(z) = \Big(\frac{2z}{|z|^2 + 1}, \frac{|z|^2 -1}{|z|^2 +1} \Big).$$
\noindent In the following, we will systematically identify $\hat{\C}$ and $\S^2$ {\it via} this  stereographic projection (justifying the term {\it Riemann sphere} for $\hat{\C}$.
\end{Myitemize}

\noindent The three surfaces $\C$,  ${\hat \C}$ and $\H$ admit a canonical structure of simply connected Riemann surface. They are the only possible ones, as claimed by the famous following result~:

\begin{theorem} {\bf -~Riemann uniformization Theorem~-} \label{UNIF}
Every complete simply connected Riemann surface is biholomorphic to $\C$, $\H$ or $\hat{\C}$.
\end{theorem}

\noindent When no confusion is possible, we will denote each of these spaces by the generic letter $\bf H$.

\subsection{Fundamental results on compact Riemann surfaces}
 We  mention here three fundamental results in the theory of Riemann surfaces (without proof). Two of them will be useful for the rest of these notes.

\begin{theorem} \label{RAC}
Every  (compact oriented) Riemann surface $X$ admits
 non constant meromorphic maps
$f: X \to \C$.
\end{theorem}

\begin{corollary}\label{RAC2}
Every  (compact oriented) Riemann surface $X$ admits
 non constant holomorphic maps into $\hat{\C}$
$$f: X \to \hat{\C}.$$
\end{corollary}

\begin{corollary}\label{RAC21}
Every  (compact oriented) Riemann surface $X$ admits
a (generally ramified) holomorphic covering over $\hat{\C}$.
\end{corollary}

\begin{theorem}\label{RAC3}
Let $X$ be a (compact oriented) Riemann surface. Then, there exists an irreducible polynomial $P \in \C[X,Y]$ such that $X$ is isomorphic to the compactification of the regular points  of the algebraic curve of equation $P(x,y)=0$.
\end{theorem}

\begin{theorem}\label{RAC31}
Let $X$ be a (compact oriented) Riemann surface. Then there exists an holomorphic embedding of $X$ into ${\C}P(3)$.
\end{theorem}

\noindent Although we will not use Theorem \ref{RAC31} in the rest of these notes, we remark that it implies that any (compact connected) Riemann surface  appears as a  $2$-dimensional real  surface  minimally embedded in the projective space (of real dimension $6$)  ${\C}P(3)$.

\section{The structure of  $\mbox{Aut}(\bf H)$}
The following result describes the group $\mbox{Aut}(\bf H)$  of (biholomorphic) automorphisms of $\bf H$~:

\begin{theorem} \label{UNIR}
\begin{enumerate}
\item \label{UNIR1}
One has~:
$$\mbox{\rm Aut}(\C)   = \{z \to az+b, a, b \in \C \};$$

\item \label{UNIR2}
One has~:
$$
\begin{aligned}
\mbox{\rm Aut}(\H) &= \{z \to \frac{az+b}{cz+d}, a,b,c,d \in \R, ad-bc=1\} \\
                   & \simeq \{ \begin{pmatrix}
                      a & b \\
                      c & d
                      \end{pmatrix}
                      ;a,b,c,d \in \R, ad-bc=1\}\\
                     & = PSL(2,\R).
\end{aligned}
$$

\item \label{UNIR3}
One has~:
$$
\begin{aligned}
\mbox{\rm Aut}(\hat{\C}) &= \{z \to \frac{az+b}{cz+d}, a,b,c,d \in \C, ad-bc=1\} \\
                         & \simeq \{ \begin{pmatrix}
                      a & b \\
                      c & d
                      \end{pmatrix}
                      ;a,b,c,d \in \C, ad-bc=1\}\\
                         &= PSL(2,\C)\simeq PGL(2,\C).
\end{aligned}
$$
\end{enumerate}
\end{theorem}

\noindent An element of $\mbox{\rm Aut}(\C)$ is called an {\it affine map} (or an {\it homothecy} if $b=0$), an element of $\mbox{\rm Aut}(\H)$ or $\mbox{\rm Aut}(\hat{\C})$ is called an {\it homography}. An element of $\mbox{\rm Aut}(\hat{\C})$ is also called a {\it Moebius transformation}. For further use, we give the following crucial result.

\begin{theorem} \label{HOMOG}
Let $a$, $b$, $c$ be three (distinct) points of $\hat{\C}$. Then there exists a unique  Moebius transformation sending $0$ to $a$, $1$ to $b$ and $\infty$ to $c$.
\end{theorem}

%
%

\noindent Now, one can build an action of the group $PSL(2,\C)$ on $\hat{\C}$ as follows~: for any matrix
$\begin{pmatrix}
a & b \\c & d
\end{pmatrix}
$, where
 $(a,b,c,d \in \C, ad-bc=1)$,

$$\begin{pmatrix}
a & b \\c & d
\end{pmatrix}.z=\frac{az+b}{cz+d}.
$$

\noindent  In the next sections, we will systematically use Proposition \ref{ISOM1} and Theorem \ref{UNIR}~:
${\rm Aut}(\hat{\C})\simeq PGL(2,\C)$.

\section{Classification of Riemann surfaces}
\subsection{Description of Riemann surfaces with respect to their genus}
 If $X$ is any Riemann surface,
$X$ is the quotient of its universal covering $\bf H$ by a subgroup of $\mbox{Aut}({\bf H})$ acting
 freely an discontinuously on ${\bf H}$. Since the covering is holomorphic, the covering automorphisms are holomorphic and then,  Riemann surfaces can be classified as follows~:

\begin{theorem} \label{UNIF2}
Let $X$ be a  Riemann surface of genus ${\bf g}_X$.
\begin{enumerate}
\item \label{UNIF3}
If  ${\bf g}_X=0$ then   $X=\hat{\C}$.
\item \label{UNIF4}
If ${\bf g}_X=1$, then $X$  is a quotient
 $\C \slash \Gamma$, where $\Gamma$ is a lattice $u\Z \oplus v\Z$, acting on $\C$ by translations, where $u$ is a complex number, $v$ is a nonzero complex number, such that $\frac{u}{v}$ is not a real number.
\item \label{UNIF5}
If ${\bf g}_X >1$, then $X$ is a quotient
 $\H \slash \Gamma$, where $\Gamma$ is a Fuchsian subgroup of $PSL(2,\R)$ (that is, a subgroup acting freely and properly discontinuously on $\H$).
\end{enumerate}
\end{theorem}

\subsection{Modular spaces}
 Let $X$ be a closed oriented surface, and ${\cal M}$ be the set of conformal structures ${\cal C}$ on $X$. Let $\sim$ be the equivalence relation defined on ${\cal M}$ as follows~: Two conformal structures ${\cal C}_1$ and ${\cal C}_2$ on $X$ are equivalent when there exists a conformal diffeomorphism
 $$\phi :(X, {\cal C}_1) \to (X, {\cal C}_2).$$

 \begin{definition}
 The quotient space ${\cal M}/\sim$ is called the {\rm modular space} of conformal structures of $X$.
 \end{definition}

 \noindent The modular spaces of conformal structures of a given genus have a structure of complex manifold. More precisely, the modular space $\mathcal{M}_{{\bf g}}$ (of conformal structures defined on a -~closed oriented~- surface of genus ${\bf g} >1$) has a structure of a complex manifold of dimension $3g-3$, and the modular space $\mathcal{M}_{1}$ can be identified to $\H/PSL(2,\Z)$. For ${\bf g}=0$, one has the following result~:

\begin{theorem} \label{MOD}
The modular space $\mathcal{M}_{0}$ is reduced to a point.
\end{theorem}

 \noindent Although we don't give a direct proof of Theorem \ref{MOD}, we remark that it is an easy consequence of the Riemann-Roch theorem~: On any (oriented closed) surface $X_0$ of genus $0$, there exists  a meromorphic function with one pole of degree one.

  \vspace{10pt}

 \noindent  On the other hand, Theorem \ref{MOD} means that $X_0$ is conformally equivalent to the Riemann sphere $\tilde{\C}$. Consequently, we can call a (closed oriented) Riemann surface  $X_0$  with genus $0$, {\it the} Riemann sphere. Concretely,  Theorem \ref{MOD} means that if ${\cal C}_1$ and ${\cal C}_2$ are two conformal structures on $X_0$, there exists a conformal diffeomorphism (or a biholomorphism preserving the orientation) from $(X_0, {\cal C}_1)$ to $(X_0, {\cal C}_2)$.

 \section{Riemannian surfaces}
  A Riemannian surface $(X,g)$ is a $2$-dimensional real (smooth) surface, endowed with a (smooth) Riemannian metric $g$. By definition, this means that $X$ is endowed with an atlas ${\cal A}=\{(U_{\alpha}, \phi_{\alpha}, V_{\alpha}\subset \R^2)\}$ where each $V_{\alpha}$ is endowed with a metric (symmetric positive definite bilinear form) $g_{\alpha}$ such that
 the transition functions
 $$\phi_{\beta} \circ \phi_{\alpha}^{-1} : (\phi_{\alpha}(U_{\alpha} \cap U_{\beta}), g_{\alpha}) \to (\phi_{\beta}(U_{\alpha} \cap U_{\beta}, g_{\beta}))$$
 are isometries.
 In the following, as usual, we identify $g$ on $U_{\alpha}$ and its local representation $g_{\alpha}$ on $V_{\alpha}$.

\subsection{Isothermal coordinates}
In local coordinates $(x,y)$ in each $U_i$, the metric $g$ defined on a Riemann surface can be written as follows~:
$$g= Edx^2 + 2Fdxdy + Gdy^2,$$
where $E$, $F$, $G$ are real valued functions of the variables $x$ and $y$.
Using complex coordinates $z=x+iy$,
$$g=\alpha|dz + \mu d\bar{z}|^2,$$
where $\alpha$, $\mu$, are real valued functions of the variables $z$ and $\bar{z}$.
Gauss (in the analytic case), Korn and Lichtenstein (in the smooth case) proved that it is always possible to find a (local) system of coordinates $(u,v)$ on $U_i$ so that
\begin{equation}\label{ISOTT}
g= e^{\lambda}(du^2 + dv^2),
\end{equation}
or, using complex coordinates $w=u+iv$,
$$g= e^{\lambda}dw d\bar{w}.$$
Such coordinates are called {\it isothermal coordinates}. This result can be stated as follows~:

\begin{theorem}
Let $(X,g)$ be a Riemannian surface. Then, around each point $p \in X$, there exists a chart $(U,\phi, V)$ and a smooth function $\lambda$ such that $g|U=e^{\lambda}g_{\E^2}$,
where $g_{\E^2}$ denotes the standard scalar product on $V \subset \E^2$.
\end{theorem}

\noindent A chart satisfying \ref{ISOTT} is called an {\it isothermal chart}.

\vspace{10pt}

\noindent The following result is obvious but important~:
\begin{proposition} \label{PROPCONF}
 The transition functions of the atlas ${\cal A}$  defined on an (oriented) Riemannian surface $(X,g)$, restricted to isothermal charts are conformal maps
 $$\phi_{\beta} \circ \phi_{\alpha}^{-1} :  \phi_{\alpha}((U_{\alpha} \cap U_{\beta}, g_{\E^2})) \to \phi_{\beta}((V_{\alpha} \cap V_{\beta}, g_{\E^2})).$$
\end{proposition}

\noindent {\footnotesize {\bf Proof of Proposition \ref{PROPCONF}~-}
Indeed, let $\phi = \phi_{\beta}^{-1} \circ \phi_{\alpha}$. Denoting for all $\alpha \in I$, $g_{\alpha} = g|U_{\alpha}$, we have
$$g_{\alpha}=e^{\lambda_{\alpha}}g_{\E^2}, g_{\beta}=e^{\lambda_{\beta}}g_{\E^2},$$
where $\lambda_{\alpha}$ and $\lambda_{\beta}$ are smooth functions, because $(U_{\alpha}, \phi_{\alpha})$ and $(U_{\beta}, \phi_{\beta})$ are isothermal charts. On the other hand, by definition of $g$,
$$\phi^*(g_{\beta}) =g_{\alpha}.$$
We deduce that
$$e^{\lambda_{\alpha}}g_{\E^2}
=g_{\alpha}=\phi^*(g_{\beta})=\phi^*(e^{\lambda_{\beta}}g_{\E^2})=e^{\lambda_j}\phi^*(g_{\E^2}),$$
from which we deduce that
$$\phi^*(g_{\E^2})=e^{\lambda_{\alpha} - \lambda_{\beta}}g_{\E^2},$$
implying that the transition functions $\phi = \phi_{\beta}^{-1} \circ \phi_{\alpha}$ are conformal.
}

    \subsection{Conformal class of a metric}
    Two Riemannian metrics $g$ and $h$ defined on a  surface $X$ are called {\it conformal} if
    $g = e^{\lambda}h$,
where $\lambda : X \to \R$ is a $C^{\infty}$ function. One can classify the Riemannian metrics by defining an equivalence relation as follows: Two Riemannian metrics $g$ and $h$ defined on $X$ are equivalent if they are conformal. For further use, we state the following lemma that gives the relation between the (Gaussian) curvatures of two conformal metrics.

\begin{lemma} \label{LEMMACURVCONF}
 Let $(X,g)$ be a closed oriented Riemannian surface with curvature $k_g$, and $h$ a metric on $X$ conformal to $g$ with (Gaussian) curvature $k_h$~:  There exists a smooth function $\lambda$ on $X$ such that
$$h=e^{\lambda}g.$$
Then,
\begin{equation} \label{COURBURE}
k_{h}=e^{-\lambda}(k_g- \Delta_g \lambda),
\end{equation}
where $\Delta_g$ is the Laplacian of $g$.
\end{lemma}

\subsection{Riemannian surfaces {\it versus} Riemann surfaces}
The link between Riemannian surfaces and Riemann surfaces can be summarized as follows~:

\begin{theorem} \label{COCO5}
Let $X$ be an (differentiable) surface. It is equivalent to endow $X$ with a complex structure or to endow it with an orientation and a conformal equivalence class of Riemannian metrics.
\end{theorem}

\noindent By a {\it complex structure}, we mean a (maximal) atlas whose transition functions are conformal.

\vspace{10pt}

\noindent {\footnotesize {\bf Sketch of Proof of Theorem \ref{COCO5}~-}
 Let us describe now the main steps of the proof of Theorem \ref{COCO5} and how the bijection is built.
\begin{itemize}
\item
 Let us show how a conformal structure on $X$ determines a natural conformal class of Riemannian metrics on $X$~: On an atlas $\mathcal{A}$, one defines  a locally finite  partition of unity $\mu_{\alpha}$, the support of each $\mu_{\alpha}$ being included in a chart $(U_{\alpha}, \phi_{\alpha}, V_{\alpha})$ of $\mathcal{A}$. One endows each $U_{\alpha}$ with the natural Euclidean metric $\phi_{\alpha}^{-1}(g_{\alpha})$, where $g_{\alpha}$ is the restriction
  on $V_{\alpha}$ of the Euclidean metric on  $\C \simeq \E^2$. Then, we define on $X$
the Riemannian metric $g=\sum \mu_{\alpha}\phi_{\alpha}^{-1}(g_{\alpha})$ on $X$. Any metric $h =e^{\lambda}g$ conformal to $g$ (where $\lambda$ is a smooth function) can be obtained by the same process (multiplying each $\mu_i$ by $e^{\lambda}$). Consequently, we have associated to any Riemann structure on $X$ a conformal class of Riemannian metrics. Moreover,
 one proves that the class of  Riemannian metrics giving rise to a given  complex structure on $X$ is exactly a conformal class of Riemannian metrics.

\item
Conversely, if $(X,g)$ is an (oriented) Riemannian surface, one  builds an isothermal  atlas $\mathcal{A} =\{(U_{\alpha}, \phi_{\alpha}, V_{\alpha}), \alpha \in I\}$ of $X$. Then we apply Proposition \ref{PROPCONF}. We conclude immediately that $X$ admits a structure of Riemann surface. By construction, two conformal metrics give rise to the same conformal structure.
\end{itemize}

\noindent The previous correspondences are inverse one to each other, Theorem \ref{COCO5} is proved.
}

\begin{definition}
A metric $g$ defined on a Riemann surface $X$ is said to be {\rm compatible} with its conformal structure if this conformal structure is induced by $g$.
\end{definition}

\noindent We conclude this section by the following remark~: Let $g_1$ ({\it resp.} $g_2$) be a Riemannian metric on a (closed oriented) surface  $X_{0}$ of genus $0$. Let ${\cal C}_1$ ({\it resp.} ${\cal C}_2$) be the conformal structure associated to  $g_1$ ({\it resp.} $g_2$). From  Theorem \ref{MOD}, we know that there exists a conformal diffeomorphism
 $$\phi : (X_0,g_1) \to (X_{0},g_2).$$

 \noindent Using Theorem \ref{MOD} and Theorem \ref{COCO5}, we get  the following Riemannian result~:

\begin{corollary} \label{CORCONF}
Let $g_1$ ({\rm resp.} $g_2$) be a Riemannian metric on a (closed oriented) surface  $X_{0}$ of genus $0$. Let ${\cal C}_1$ ({\rm resp.} ${\cal C}_2$) be the conformal structure associated to  $g_1$ ({\rm resp.} $g_2$). Then, there exists a  biholomorphism
 $$\phi : (X_0,g_1) \to (X_{0},g_2),$$
  that is, there exists a  diffeomorphism $\phi$ of $X_0$ and a smooth function $\lambda$ on $X_{0}$ such that, $\phi^*(g_2)=e^{\lambda}g_1$.
\end{corollary}

\noindent Indeed, $g_1$ ({\it resp.} $g_2$) induces a conformal structure ${\cal C}_1$ ({\it resp.} ${\cal C}_2$). Since Theorem \ref{MOD} claims that the modular space of $X_0$ is reduced to a point, there exists a conformal diffeomorphism
 $$\phi : (X_0, {\cal C}_1) \to (X_{0}, {\cal C}_2),$$
 and then, there exists a smooth function $\lambda$ such that $\phi^*(g_2)=e^{\lambda}g_1$.

\section{Metrics of constant Gaussian curvature} \label{KJOI77}
\subsection{The simply connected case} \label{KJOI9}
When $X$ is simply connected, one can give an explicit description of the metrics with constant curvature on $X$. We begin with the famous theorem of Cartan~:

\begin{theorem} {\bf -~The uniformization theorem of Cartan in Riemannian geometry~-} \label{THCAR}
Let $(X,g)$ be a connected, complete, simply connected, oriented Riemannian surface with constant Gaussian curvature $k$.
\begin{itemize}
\item
If $k=1$, then $(X,g)$ is isometric to the round sphere  of radius $1$.
\item
If $k=0$, then $(X,g)$ is isometric to the Euclidean plane.
\item
If $k=-1$, then $(X,g)$ is isometric to the hyperbolic plane (endowed with its canonical metric of constant curvature $-1$).
\end{itemize}
\end{theorem}

\noindent Remark that these isometries are not unique since they are parametrised by isometric automorphisms of $(X,g)$. On the other hand, the three spaces described in Theorem \ref{THCAR} (the sphere, the plane, the hyperbolic plane), admits by Theorem \ref{UNIF}  a canonical complex structure. We now describe these three situations.
\begin{enumerate}
\item
Let us study the case $k=0$, and identify the plane with $\C$. The metric $g_{\C}=|dz^2|$ is the canonical flat Riemannian metric on $\C$. Since any biholomorphism of $\C$ is an affine map $$f: z \to az +b, a \in \C, b \in \C,$$
we can write $f^{*}g_{\C}=e^{\lambda} g_{\C}, \lambda \in \R$, from which we deduce a family of flat  Riemannian metrics on $\C$ associated to the canonical complex structure, conformal to $g_0$, parametrised by $\R^{*}_{+}$. So, up to a scaling constant, there exists a canonical  flat Riemannian metric on $\C$ associated to its canonical complex structure.

\item
Let us study the case $k=-1$, identifying the hyperbolic plane with  $\H$. The Riemannian metric
$$g_{\H}=\frac{|dz|^2}{[\mbox{{\rm im}}(z)]^2}$$  has constant Gaussian curvature $-1$. The crucial remark is that $g_{\H}$ is invariant
by $\mbox{\rm Aut}(\H)$. In other words, any biholomorphism of $\H$ is a $g_{\H}$-isometry. We call  $g_{\H}$ the {\it canonical metric} of constant curvature $-1$ on $\H$.

\item
Let us study the case $k=1$. {\it Via}  the stereographic projection $s$, we identify the sphere $\S^2$ as before with $\hat{\C}$. The sphere $\S^2$ admits a Riemannian metric $g_{\S^2}$ of constant curvature  $1$, induced by the standard scalar product of $\E^3 \simeq \C \times \R$~:
$$g_{\S^2}=\frac{4|dz|^2}{(1 + |z|^2)^2}.$$
On $\hat{\C}$,  a simple computation gives
\begin{equation} \label{CANME}
g_{\hat{\C}}= \frac{4|du|^2}{(1+|u|^2)^2},
\end{equation}
where $u=\frac{1}{z}$.
However, $g_{\hat{\C}}$ is not invariant by $\mbox{\rm Aut}(\hat{\C})$.  Consequently $g_{\hat{\C}}$ is not characterised by the Riemann structure of $\hat{\C}$. But the curvature of $g_{\hat{\C}}$ is preserved by $\mbox{\rm Aut}(\hat{\C})$. More precisely, we have the following lemma~:

 \begin{lemma} \label{ISOMH}
 \begin{enumerate}
 \item
 If $h \in \mbox{\rm Aut}(\hat{\C})$, $h^*(g_{\hat{\C}})$ is a Riemannian metric with constant curvature $1$.
 \item \label{ITTEE}
The family of Riemannian metrics $\tilde{g}$ with constant Gaussian curvature $1$ conformal to the canonical metric $g_{\hat{\C}}$  on $\hat{\C}$,  is naturally endowed with a structure of homogenous space isomorphic to   $SO(3)\backslash \mbox{\rm Aut}(\hat{\C}) \simeq SO(3)\backslash PSL(2,\C)$.

\item
More generally, if $g$ is any Riemannian metric on $\hat{\C}$, the family of Riemannian metrics $\tilde{g}$ with constant Gaussian curvature $1$ conformal to $g$, is naturally endowed with a structure of homogenous space isomorphic to   $SO(3) \backslash \mbox{\rm Aut}(\hat{\C}) \simeq SO(3) \backslash PSL(2,\C)$.
  \end{enumerate}
 \end{lemma}

 \noindent {\footnotesize {\bf Proof of Lemma \ref{ISOMH}~-}
 \begin{enumerate}
 \item
 If $h \in \mbox{\rm Aut}(\hat{\C})$, $h^*(g_{\hat{\C}})$ is still a metric with constant curvature $1$, but generally, $h^*(g_{\hat{\C}}) \neq g_{\hat{\C}}$. In other words, any biholomorphism of $\hat{\C}$ preserves the (constant) curvature of $g_{\hat{\C}}$ but does not preserve $g_{\hat{\C}}$ in general.

 \item
 Let $\tilde{g}$ be any metric of constant curvature $1$ on $\hat{\C}$. Then, by Theorem \ref{THCAR}, there exists an isometry
 $$\psi: (\hat{\C},  \tilde{g})  \to (\hat{\C}, g_{\hat{\C}}).$$
 On the other hand, if we suppose that $\tilde{g}$ is in the conformal class of $g_{\hat{\C}}$,  there exists a smooth function $\lambda$ such that
  $\tilde{g}=e^{\lambda}g_{\hat{\C}}$. We deduce that
  $$\psi^{*}g_{\hat{\C}}=e^{\lambda}g_{\hat{\C}},$$
 that is, $\psi$ is a conformal map from  $(\hat{\C},g_{\hat{\C}})$ to $(\hat{\C},g_{\hat{\C}})$~: $\psi \in \mbox{Aut}(\hat{\C}) \simeq PSL(2,\C)$. Now, let $\varphi \in \mbox{Aut}(\hat{\C})$. We have
 $$(\varphi \circ \psi)^{*}g_{\hat{\C}}=(\psi^* \circ \varphi^*)g_{\hat{\C}}=\psi^*g_{\hat{\C}},$$ for all  $\psi \in \mbox{Aut}(\hat{\C})$, if and only if $\varphi \in  SO(3)$. The conclusion follows.

 \item The last item is an easy generalization of item \ref{ITTEE}~: If $g$ is any metric on $\hat{\C}$, Theorem \ref{PKK} claims that there exists a metric $\tilde{g}$ of constant curvature $1$ conformal to $g$. If $\tilde{g}'$ is another  metric of constant curvature $1$ conformal to $g$, then $\tilde{g}'$ is conformal to $\tilde{g}$. We apply item \ref{ITTEE}, and the conclusion follows.

\end{enumerate}
}
\end{enumerate}

\subsection{A general uniformization theorem in Riemannian geometry} \label{PKK2}
 We consider now closed Riemannian surfaces with any genus. The {\it Poincare-Klein-Koebe uniformization theorem} claims that in each conformal class of metrics on a surface $X$, there exists a metric with constant Gaussian curvature $1$, $0$ or $-1$~:

\begin{theorem} {\bf -~Poincare-Klein-Koebe uniformization theorem~-} \label{PKK}
Let $(X,g)$ be a  closed oriented Riemannian surface of genus ${\bf g}_X$. Then, there exists a Riemannian metric $\tilde{g}$ conformal to $g$ with constant Gaussian curvature $-1$, $0$, or $1$.
\begin{enumerate}
 \item
If ${\bf g}_X > 1$,  the Gaussian curvature is $-1$, and $g$ is unique.
\item
If ${\bf g}_X =1$,  the Gaussian curvature is $0$, and $g$ is unique up to a scaling constant.
\item \label{TROIS}
If ${\bf g}_X=0$, there exists a family of Riemannian metrics $\tilde{g}$ conformal to $g$ with constant Gaussian curvature  $1$.
\end{enumerate}
\end{theorem}

\noindent We will improve Theorem \ref{PKK} item \ref{TROIS} in section \ref{KJOI9} Lemma
 \ref{ISOMH} by describing the geometry of the set of metrics  with constant curvature $1$ that are conformal to $g$.
The proof of  Theorem \ref{PKK} is based on Lemma \ref{LEMMACURVCONF}.
To find a metric of constant curvature conformal to $g$, one solves equation  \ref{COURBURE} with $k_{g'}=-1,0$ or $1$, that is, one looks for a smooth function $\lambda$ satisfying \ref{COURBURE}.
\vspace{10pt}

\noindent {\footnotesize {\bf Sketch of proof of  Theorem \ref{PKK}~-}
We only give here indications of the proof in the simplest case of genus $1$.
 Let $X_1$ be a surface of genus $1$ endowed with a metric $g$ of curvature $k_g$. Let $g'=e^{\lambda}g$, where $\lambda$ is a smooth function on  $X_0$. From equation \ref{COURBURE}, we deduce that $k_{g'}=0$ if and only if
 \begin{equation} \label{COURBURE2}
\Delta_g \lambda=k_g,
\end{equation}
 By using the classical  theory of autoadjoint operators, we solve equation \ref{COURBURE2}~: Up to a constant, it admits a unique solution.
}


%
\subsection{Metrics of constant curvature on a Riemann surface of any genus}
 If we {\it a priori} deal with a closed Riemann surface (of any genus),  the results of subsection \ref{PKK2} can be rephrased as follows~:

\begin{theorem} \label{THFI}
Let $X$ be a (closed oriented) Riemann surface  of genus ${\bf g}_X$. Then,
\begin{itemize}
 \item
If ${\bf g}_X > 1$,   there exists a {\rm unique} Riemannian metric $\tilde{g}$  with constant Gaussian curvature $-1$, compatible with the complex structure and invariant by $\mbox{\rm Aut}(X)$.
\item
If ${\bf g}_X =1$,   there exists a  Riemannian metric $\tilde{g}$ with constant Gaussian curvature $0$, {\rm unique up to a scaling constant}, compatible with the complex structure and invariant by $\mbox{\rm Aut}(X)$.
\item \label{THFI3}
If ${\bf g}_X =0$, and $g$ is a Riemannian metric on $X$,  the family of  Riemannian metrics $\tilde{g}$ with constant Gaussian curvature $1$ conformal to $g$ is naturally endowed with a structure of homogenous space isomorphic to   $\mbox{\rm Aut}(SO(3) \backslash \hat{\C}) \simeq SO(3) \backslash PSL(2,\C)$.
\end{itemize}
\end{theorem}

\noindent We deduce from Theorem \ref{THFI} that there exists  a natural homogenous space of metrics with constant Gaussian curvature $1$ on $\hat{\C}$~: We consider the metric $g_{\S^2}$ of constant Gaussian curvature $1$ on the round sphere $\S^2$ of radius $1$. We identify $\S^2$ with $\hat{\C}$ by help of the stereographic projection, and consider the metric $g_{\hat{\C}}$ on  $\hat{\C}$ deduced from $g_{\S^2}$ by this identification.  We then consider the  homogenous space of  Riemannian metrics $\tilde{g}$ with constant Gaussian curvature $1$ conformal to $g_{\hat{\C}}$.

\subsection{Metric on $\hat{\C}$ associated to a finite subgroup of $\mbox{\rm Aut}(\hat{\C})$}
The following proposition builds a canonical Riemannian metric on $\hat{\C}$ associated to any finite subgroup $K$ of ${\rm Aut}(\hat{\C})$  (by Theorem \ref{SGSO}, we know that $K$  is conjugate to a (finite) subgroup of $SO(3)$ canonically embedded in $\mbox{\rm Aut}(\hat{\C})$ by Proposition \ref{SO3}).

 \begin{proposition}\label{SGSO21}
 Let $K$ be a finite subgroup of ${\rm Aut}(\hat{\C})$. Then,
   $\hat{\C}$ admits a canonical Riemannian metric $\tilde{g}_K$ invariant by $K$, that coincides with the metric $g_{\hat{\C}}$  if $K$ is a subgroup of $SO(3)$.
 \end{proposition}

 \noindent {\footnotesize {\bf Proof of Proposition \ref{SGSO21}~-}
We still denote by  $g_{\hat{\C}}$ the canonical metric on $\hat{\C}$ defined by \ref{CANME}. We define $\tilde{g}_K$ as follows~:
For all $u,v$ in $T\hat{\C}$,
\begin{equation}\label{KAAP}
\tilde{g}_K(u,v)= \frac{1}{{\rm card}(K)}\sum_{h \in K}g_{\hat{\C}}(dh(u), dh(v)).
\end{equation}
For each $h \in K$, Lemma \ref{ISOMH} implies that the map
$$(u,v) \to g_{\hat{\C}}(dh(u), dh(v))$$
is a metric of constant curvature $1$.
Then, $\tilde{g}$ is a metric clearly invariant by $K$. The rest of Proposition \ref{SGSO21} is clear.
}


\section{Canonical Riemann structure on a  polyhedron} \label{CANRI}
By an (abstract) oriented Euclidean polyhedron, we mean an oriented topological surface obtained as the union of a finite set of disjoint polygonal domains of the Euclidean plane $\E^2$, after the identification of some of their edges and vertices. Such polyhedra are piecewise linear and admit on each of their face a canonical (Euclidean) flat metric with potential singularities at the vertices.
In this section we will build a canonical Riemann structure on any (abstract) Euclidean polyhedron. More precisely, we can claim~:

\begin{theorem} \label{PLT}
Any Euclidean polyhedron admits a canonical conformal structure.
\end{theorem}

\noindent {\footnotesize {\bf Proof of Theorem \ref{PLT}~-}
 Let us build a conformal atlas on the polyhedron $P$~:

 \begin{Myitemize}
 \item
 First of all, we consider any point $p$ belonging to the interior of a face $f$ or to the interior of an edge $e$ adjacent to two faces $f$ and $f'$. As a domain of chart around $p$, we take any open neighborhood $U$ of $p$ in $f \cup f'$, and we send it isometrically onto a neighborhood of $0$ in the Euclidean plane $\E^2$. Denoting by $\phi$ this isometry, the triple $(U, \phi, \phi(U))$  is a chart around $p$.

 \item
 On the other hand, if $v$ is a vertex of the polyhedron, we consider a "small" neighborhood $U$ of $v$ that is the union of sectors at $p$ surrounding $p$~:  Let $e_1,\ldots, e_k$ be the sequence of edges adjacent to $v$, and $\theta_1, \ldots, \theta_k$ the angle between consecutive edges.
  \begin{itemize}
  \item
  We send isometrically the interior $\stackrel{\circ}{f_i} \cap U$ in each face $f_i$ to the interior of a sector $V_i$  in $\E^2$ whose vertex is $0$.

  \item
  Now, we apply the transformation $\delta_i: z \to z^{\alpha}$, where
  $$\alpha =\frac{2\pi}{\sum_k \theta_i}.$$
  (This transformation is well defined at every point different to $0$.)
  In such a way, the $\delta_i(\stackrel{\circ}{f_i}\cap U)$ is a sector  $V_i$ of angle $\alpha \theta_i$.

  \item
  After rotations with suitable angles, the union of the sectors  $V_i \cup e_i$ covers exactly an open neigborhood of $0$ in $\E^2$ (punctured at $0$), since the sum of the sector angles equals $2\pi$.

  \item
  Therefore, by continuity on the edges and on the vertex $p$, we have built a homeomorphism $\phi$ of $U$ onto an open neighborhood of $0$. We remark that $v$ is sent onto $0$. The triplet $U, \phi, \phi(U)$ is a chart.

  \item
  The set of charts $\{U_i, \phi_i, \phi_i(U_i), i \in I\}$ defined above define an atlas on $P$. Let us study the transition functions. Let $U_i$ and $U_j$ be two domains of charts.
    \begin{itemize}
    \item
    If $U_i$ or $U_j$  contains no vertices, the transition function
   $$\phi_j^{-1} \circ \phi_i :\phi_i(U_i \cap U_j) \to \phi_j(U_i \cap U_j)$$
   is composed of rotations, translations that are holomorphic,  and power functions $z \to z^{\alpha}$, that are  holomorphic since the origin does not belong to $\phi_i(U_i \cap U_j)$.
    \item
    If $U_i$ or $U_j$  contains a vertex, the transition function
   $$\phi_j^{-1} \circ \phi_i :\phi_i(U_i \cap U_j) \to \phi_j(U_i \cap U_j)$$
   is composed of rotations, translations that are holomorphic, and power functions $z \to z^{\alpha}$, that are bounded and  holomorphic except at $0$ (where the power function is not defined). Then, by Theorem \ref{POINT}, $\phi_j^{-1} \circ \phi_i$ can be extended to an holomorphic function on $\phi_i(U_i \cap U_j)$.
   \end{itemize}
    \end{itemize}
     \end{Myitemize}

\noindent Therefore, this construction defines a holomorphic structure on $P$.
}

\begin{figure}[ht!]
 \begin{center}
 \includegraphics[width=0.4\textwidth]{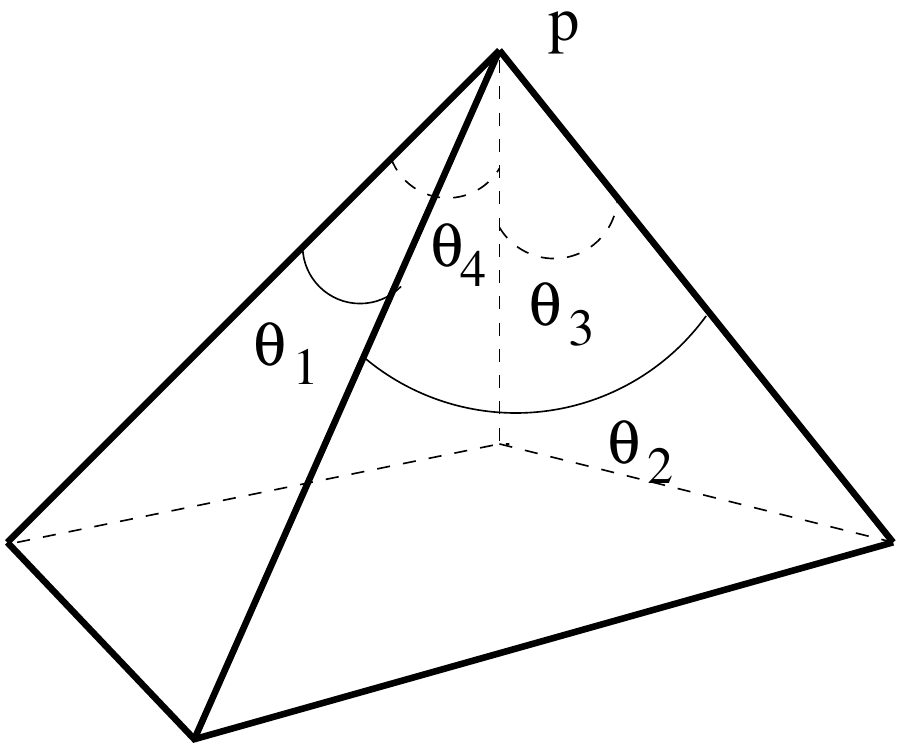}
 \end{center}
 \caption{Angles at a vertex $p$ of the triangles incident to $p$}
 \end{figure}

%
 \begin{figure}[ht!]
 \begin{center}
 \includegraphics[width=0.6\textwidth]{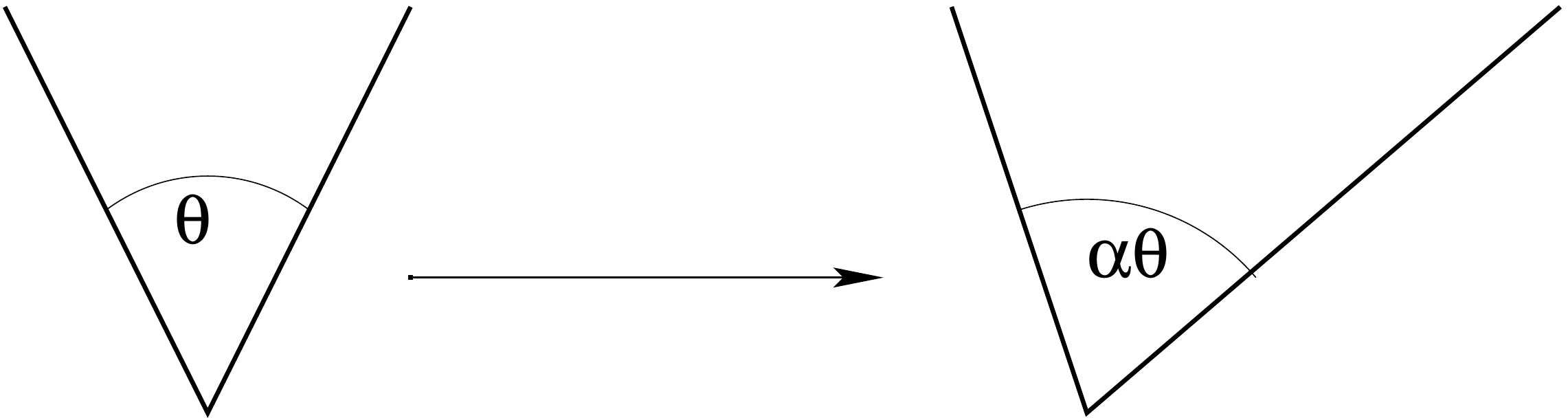}
 \end{center}
 \caption{Modifying the angle of each sector}
 \end{figure}

  \begin{figure}[ht!]
 \begin{center}
 \includegraphics[width=0.7\textwidth]{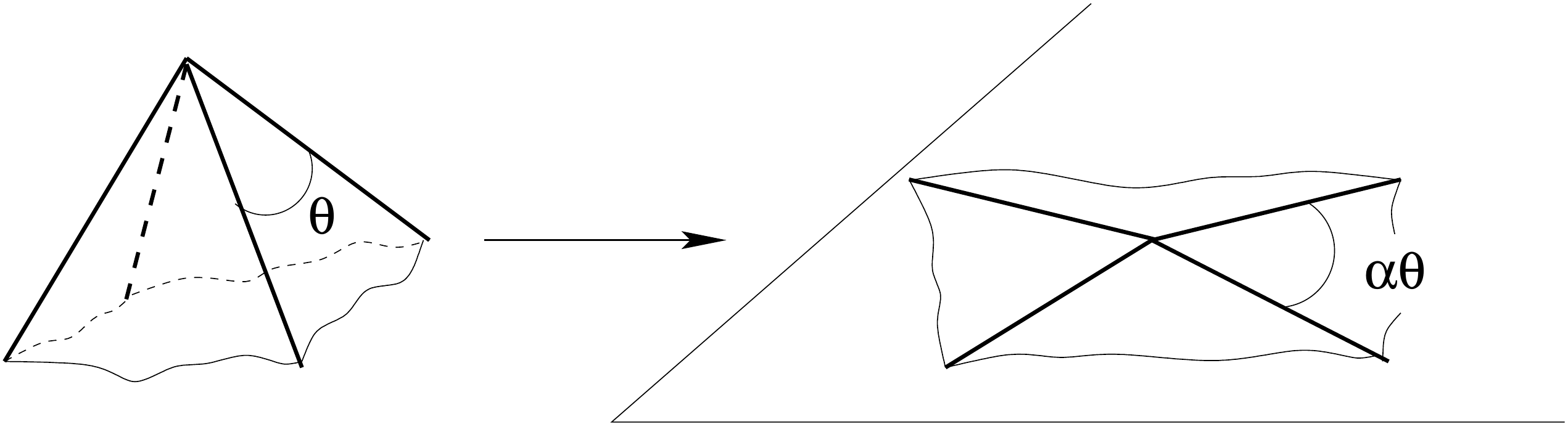}
 \end{center}
 \caption{Flattening a cone}
 \end{figure}

%
%
%

\noindent In particular, the (boundary of) any polyhedric body in $\E^3$ admits a canonical conformal structure.

\vspace{10pt}

\noindent Theorem \ref{PLT} can be extent to any triangulation on any (compact oriented) topological surface $X$. Indeed, by assigning the length $1$ to each edge of ${\cal T}$, $X$ is canonically endowed with a structure of Euclidean polyhedron whose faces are equilateral triangles (the Euclidean structure of each triangle being induced by the ones of the edges). Such a geometric structure will be called  an {\it equilateral triangulation}). We deduce~:

\begin{corollary}\label{TSD}
Any triangulation ${\cal T}$ defined on a (closed oriented) surface $X$ induces on $X$ a canonical conformal structure: The one defined by the equilateral triangulation.
\end{corollary}

\noindent The following theorem claims that the converse of Theorem \ref{PLT} is true in the following sense~:

\begin{theorem} \label{BEPOL}
Let $X$ be a (closed oriented) Riemann surface endowed with a conformal structure ${\cal C}$. Then, there exists on $X$ a Euclidean  triangulation ${\cal T}$ whose associated conformal structure is isomorphic to ${\cal C}$.
\end{theorem}

\noindent {\footnotesize {\bf Proof of Theorem \ref{BEPOL}~-}
\begin{itemize}
\item
The result is trivial if the genus of $X$ is $0$, since the modular space of $X$ is reduced to a point (see Theorem \ref{MOD}).

\item
Let $X$ be any (closed oriented) Riemann surface. We know that $X$ admits  holomorphic (generally ramified) coverings over $\hat{\C}$
$$\pi : X \to \hat{\C},$$
(see Theorem \ref{RAC21}).
We choose one of them.
Let ${\cal T}'$ be any triangulation of $\hat{\C}$ satisfying the following property~: Any singular value of $\pi$  is a vertex of ${\cal T}'$. Let us endow ${\cal T}'$ with any Euclidean  structure. By Theorem \ref{PLT}, it  induces the (unique) conformal structure of $\hat{\C}$. Moreover, the inverse image of ${\cal T}'$ by $\pi$ is a triangulation ${\cal T}$ of $X$ such that each triangle $t$ of ${\cal T}$ is in one to one correspondence with a triangle of ${\cal T}'$. Let us endow each triangle $t$ with the Euclidean metric making  the restriction of $\pi$ to  $t$  an isometry and then a conformal bijection from $t$ to  $\pi(t)$. Let us denote by $\{v_1,...,v_n\}$ the vertices of $\cal T$ and by $\{w_1,...,w_m\}$ the vertices of ${\cal T}'$. The pullback of the  conformal structure of $\hat{\C} \backslash \{w_1,...,w_m\}$ by $\pi$ is a conformal structure ${\cal C}_1$ on $X \backslash \{v_1,...,v_n\}$ isomorphic to the restriction of ${\cal C}$ on $X \backslash \{v_1,...,v_n\}$. Because the set of vertices on ${\cal T}$ is finite, ${\cal C}_1$ can be extended to $X$, and  ${\cal C}_1={\cal C}$ on $X$.
\end{itemize}
}

\section{Dessins d'enfants}\label{DESENF}
\subsection{Main definitions}
Introduction and developments on the theory of {\it dessins d'enfants} can be found in \cite{jones2016dessins} \cite{lando2013graphs} \cite{girondo2011introduction} \cite{voisin1980cartes} \cite{malgoire12cartes} \cite{grothendieck1984esquisse}. The definition of dessins may differ following the authors. For simplicity, we will use the following "standard" one:

\begin{definition} \cite{grothendieck1984esquisse} \cite{voisin1980cartes} \cite{malgoire12cartes}
A {\rm dessin d'enfants} (or simply a {\rm dessin}) is a couple $(X, {\cal D})$, where $X$ is an (oriented) topological  surface and ${\cal D}$ a  finite bicolored graph on it such that $X \backslash D$ is a finite union of disjoint topological discs.
\end{definition}

  \begin{figure}[ht!]
 \begin{center}
 \includegraphics[width=0.6\textwidth]{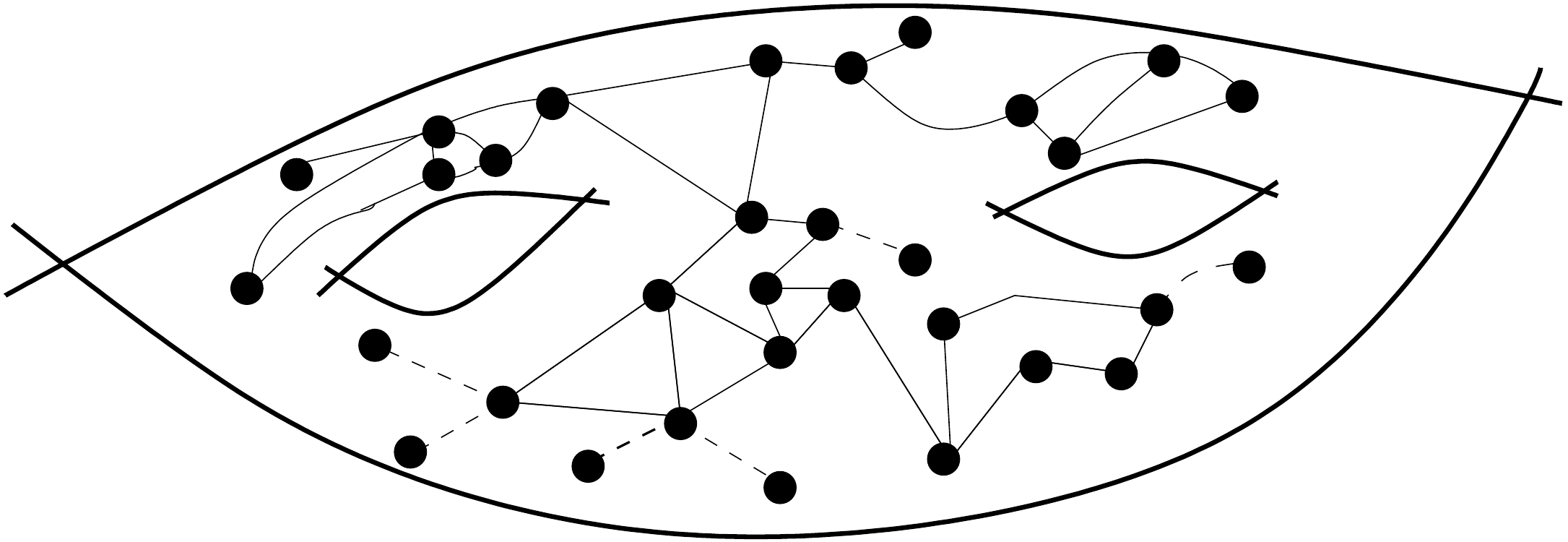}
 \end{center}
 \caption{A graph on  a surface of genus $2$}
 \end{figure}

\noindent As any finite graph, a dessin has a finite set of {\it vertices} $V$ and  a finite set of  {\it edges} $E$. As a bicolored graph,  each vertex can be colored in white or black, in such a way that the colors of two consecutive vertices on the same edge are different.
%

%

\vspace{10pt}

\noindent The definition of a dessin by A. Grothendieck is a little bit more general~: One simply considers a graph on a surface $X$ such that $X \backslash D$ is a finite union of disjoint topological discs, without the bipartite property. However, one can recover the previous definition by coloring each vertex is black, and adding a white vertex in the interior of each edge. One gets a bipartite graph whose each white vertex has valence $2$. The reader will check that many results of the following sections does not need a coloring of the dessin.

\subsection{Triangulation associated to a dessin} \label{BF1}
To each dessin $(X, {\cal D})$ is associated  a triangulation ${\cal T}$ built as follows~:
\begin{enumerate}
\item \label{PP1}
In the interior of each face $f$ of ${\cal D}$, we choose a new vertex (marked  by $*$  and called the {\it center}  of $f$).

\item \label{PP2}
From each white ({\it resp.} black) vertex $v$ of $f$, we draw a new edge joining $v$ to $*$ so that the interior of two different edges have no intersection points, and the interior of such an edge has no intersection points with the edges of $f$. This process induces a triangulation of $f$. Remark however that two adjacent triangles may have two common edges.

\item \label{PP3}
By continuing this process for each face $f$ of ${\cal D}$, we get a triangulation of $X$. Since $X$ is oriented, we can color any triangle in white and its adjacent ones in black, getting two classes of triangles of ${\cal T}$, (of types $+$ and type $-$ for instance).
\end{enumerate}

  \begin{figure}[ht!]
 \begin{center}
 \includegraphics[width=0.5\textwidth]{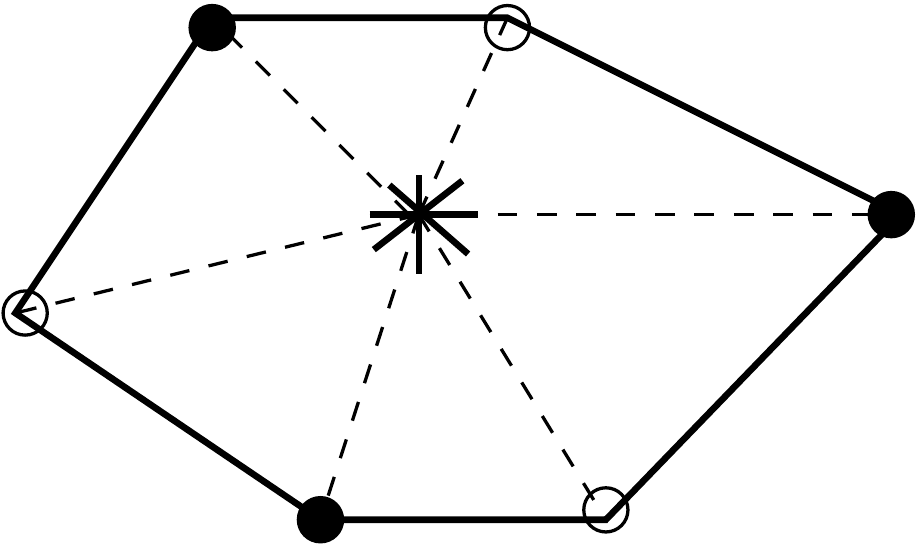}
 \end{center}
 \caption{Triangulation of a bicolored face}
 \end{figure}

  \begin{figure}[ht!]
 \begin{center}
 \includegraphics[width=0.3\textwidth]{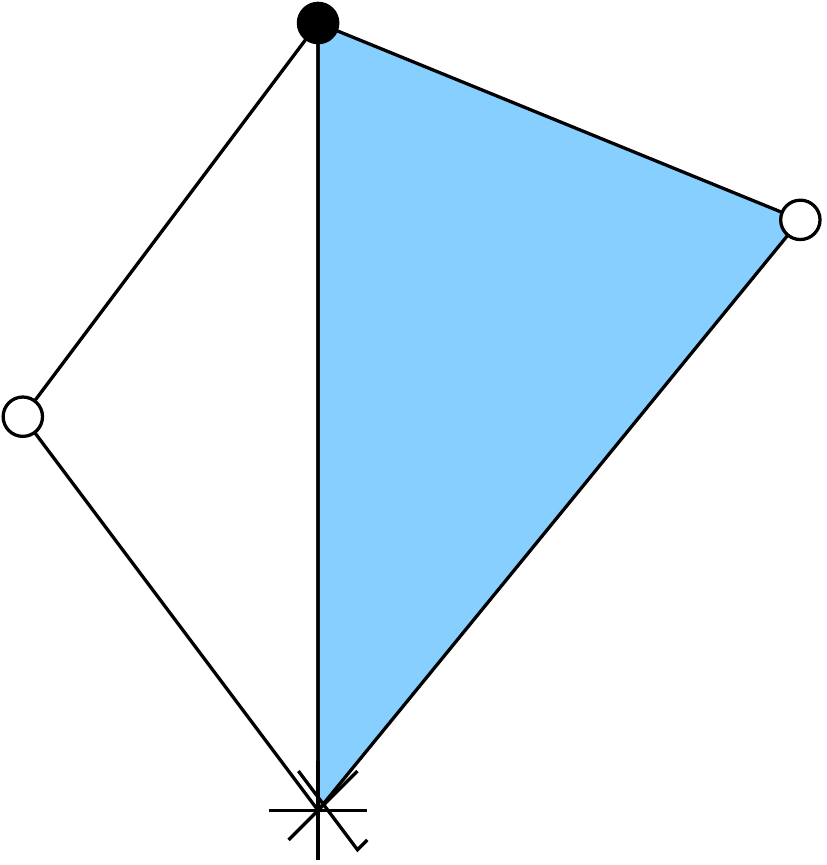}
 \end{center}
 \caption{A butterfly}
 \end{figure}

  \begin{figure}[ht!]
 \begin{center}
 \includegraphics[width=0.5\textwidth]{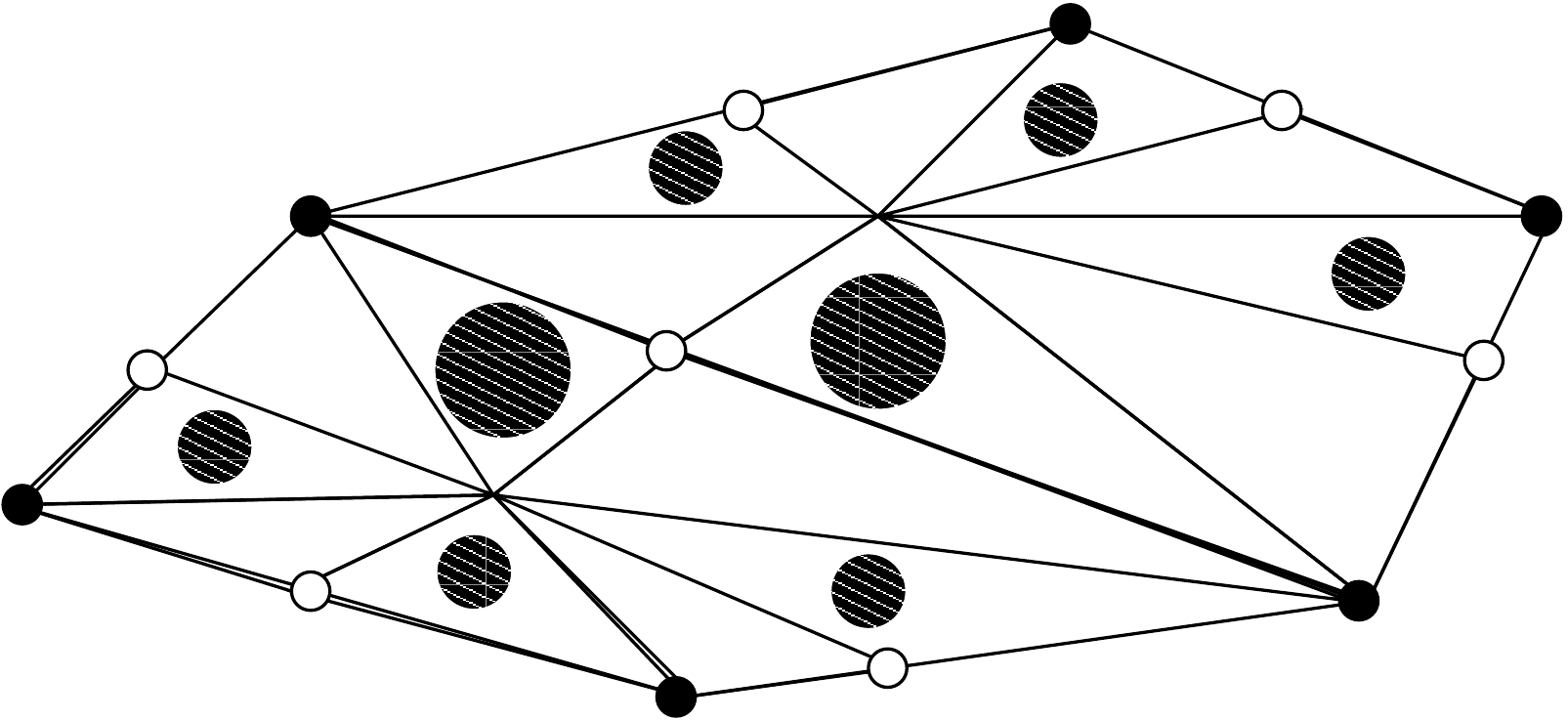}
 \end{center}
 \caption{A union of butterflies}
 \end{figure}

\noindent We remark that each triangle of this new triangulation has three different vertices~: white, black, $*$; each edge is adjacent to a triangle $+$ and a triangle $-$; each face of ${\cal D}$ is the union of an even number $2p$ of triangles, with $p$ triangles of type $+$ and $p$ triangles of type $-$. Of course, this construction depends on the positions of the center of the faces $f$ and of the shape of the new edges, but we will see that it is not important for our purpose. For further use, following \cite{shabat1990drawing}, we call {\it butterfly} the union of a triangle $t^+$ and a triangle $t^-$ adjacent at an $(.-*)$-edge of the new triangulation ${\cal T}$, so that, for instance,  triangular face is the union of three butterflies.

\vspace{10pt}

\noindent Here are two examples~:
 \begin{itemize}
 \item
 Let us consider the unit $2$-sphere $\S^2$ of $\E^3$ endowed with its equator. Let $0=(0,0,1)$, $1=(0,1,0)$, $N$ be the north pole and $S$ be the south pole. Let us build the edges $0-1$  and $1-0$ on the equator. This is the simplest dessin on $\S^2$, with two faces~! Moreover, an associated triangulation is built by adding $N$ and $S$ and drawing curves from $N$ to $0$ and $1$, ({\it resp.} $S$ to $0$ and $1$).

 \item
 If ${\cal T}$ is any triangulation on a surface $X$, we color its vertices in black. We add a white vertex in the interior of each edge and a $*$-vertex in the interior of each face, we build a new triangulation by using the construction described in \ref{PP2} and \ref{PP3}.
 \end{itemize}

\section{Building  complex structures on a dessin} \label{BF2}
Our goal now is to define  explicit  complex structures on a dessin, based on the previous constructions.
\begin{itemize}
\item
In subsection \ref{SUBPO1}, we will associate a first conformal structure ${\cal C}_1$ to a topological surface $x$ endowed with a dessin ${\cal D}$.

\item
In subsection \ref{SUBPO2}, we will associate to a topological surface $x$ endowed with a dessin ${\cal D}$ a (in general ramified and not unique) covering map
$$\beta : X \to \hat{\C},$$
that induces on $X$ a  conformal structure ${\cal C}_2$.

\item
In subsection \ref{DEUXCOIN}, we will compare ${\cal C}_1$ and ${\cal C}_2$.
\end{itemize}

\subsection{Construction of ${\cal C}_1$} \label{SUBPO1}
Let $(X, {\cal D})$ be a topological (closed oriented)  surface endowed with a dessin. Let ${\cal T}$ be a triangulation obtained by the construction described in section \ref{BF1}. We endow ${\cal T}$ with the structure of  equilateral Euclidean triangulation by affecting the length $1$ to each edge of ${\cal T}$. We remark that this Riemannian structure is independent of the choice of the position of the vertex added in each face and the "shape" of the edges of ${\cal T}$, since two such Euclidean triangulations are isometric by construction. Then, we apply Theorem \ref{PLT} (or directly Corollary \ref{TSD})~: ${\cal T}$ and then $X$ is endowed with a conformal structure. We call it ${\cal C}_1$.

\subsection{Construction of ${\cal C}_2$} \label{SUBPO2}
 To build a second conformal structure ${\cal C}_2$ on a  topological (closed oriented)  surface endowed with a dessin, we follow the idea of A. Grothendieck. It needs two steps~:

 \subsubsection{Construction of a (generally) ramified covering over $\hat{\C}$}
 We build a ramified covering
  $$\beta : X \to \hat{\C}$$
   as follows~: We begin to build a triangulation ${\cal T}$ of ${\cal D}$ by defining  butterflies $t^+ \cup t^-$  (see section \ref{BF1}), in such a way that ${\cal T}$ becomes a union of butterflies. Then,
   \begin{itemize}
   \item
   one  builds  an  homeomorphism from each butterfly $b=t^+ \cup t^-$ to $\hat{\C} \simeq \S^2$, whose equator is identified with $\R \cup \infty$.
The triangle $t^+$ ({\it resp.} $t^-$) is sent homeomorphically onto the superior ({\it resp.} inferior) hemisphere by sending the boundary of $t^+$ onto the equator, such that the white vertex is sent to $0$, the black one onto $1$ and the $*$ onto  $\infty$.
We get an homeomorphism from $b$ to $S^2$.

 \item
 By building such an homeomorphism for each butterfly, we build a map
$$\beta :X \to \hat{\C} \simeq \S^2,$$
from $X$ to $\S^2$, that is locally one-one at each point of $X$ different to the vertices of ${\cal D}$~: each white vertex is sent onto $0$, each black vertex is sent onto $1$ and each $*$ is sent to $\infty$. Consequently, $\beta$ a covering from $X$ over $S^2$, depending on $\cal D$, ramified at most above the points $0$, $1$ and $\infty$. The degree of this covering is the number of butterflies, the index of ramification of $0$ is the number of white vertices, the index of ramification of $1$ is the number of black vertices, and the index of ramification of $\infty$ is the number of faces of ${\cal D}$.
\end{itemize}

\subsubsection{Construction of ${\cal C}_2$ by pullback}
We will now use a classical result on complex functions~:

\begin{proposition} \label{PIY}
Let $X_1$ be a topological surface, $X_2$ be a Riemann surface,
$$\beta : X_1 \to X_2$$
be a covering  of finite degree, ramified at a finite number of points $z_1, ..., z_k$ of $X_1$. Let $X'_1 = X_1 \backslash \{z_1, ..., z_k\}$. Then,
\begin{enumerate}
\item
$X'^1$ is canonically endowed with the complex structure defined by pulling back by $\beta$ the complex structure of $X_2 \backslash \{\beta(z_1),...,\beta(z_k)\}$.
\item
Moreover, there exists a unique Riemann structure on $X_1$ extending the one defined on $X'^1$, so that $\beta$ is holomorphic.
\end{enumerate}
\end{proposition}

\noindent Applying Proposition \ref{PIY}, with $X_1=X$, $X_2= \hat{\C}$, $\beta$ the (generally ramified) covering build in the previous paragraph, $k=3$, $\beta(z_1)=0$, $\beta(z_2)=1$, $\beta(z_3)=\infty$, we conclude  that a dessin $(X, {\cal D})$ is endowed with a  complex structure. However, it is also clear that this construction depends on $\beta$ and the choice of the vertices $*$ in each face, the shape of the edges $0-*$ and $1-*$ and on the choice {\it a priori} to send the white points on $0$ and the black ones on $1$ (we could do the converse). But, modulo an equivalence of covering and a Moebius transformation of $\hat{\C}$, this construction is unique. Moreover, by pulling back the complex structure of $\hat{\C}$  onto $X$, one endows $X$ with a structure of Riemann surface, and ${\cal T}$ is the inverse image of the triangle $0,1,\infty$ of $\hat{\C}$.

\vspace{10pt}

\noindent We know (Theorem  \ref{RAC}) that
 any (compact oriented) Riemann surface admits  a (generally ramified) covering  over ${\hat{\C}}$.
The  construction described in  Proposition  \ref{PIY} allows to build on any (compact oriented) topological surface endowed with a dessin, a structure of Riemann surface (endowed with a conformal structure ${\cal C}_2$) and a {\it particular} (generally ramified) covering over $\hat{\C}$ that is holomorphic. Such a covering ramifies at most over three points. This leads to introduce the following definition~:

\begin{definition}\label{BEY}
Let $X$ be a compact Riemann surface. A non constant meromorphic function
$$\beta : X \to \hat{\C},$$
is a {\rm Belyi function} if it ramifies at most above three points. In this case, the couple $(X,\beta)$ is called a {\rm Belyi pair}.
\end{definition}

\noindent  Usually, {\it via} a Moebius transformation of $\hat{\C}$, the three points involved in Definition \ref{BEY} can be systematically taken as $0,1,\infty$.

\vspace{10pt}

\noindent We deduce the following theorem~:

\begin{theorem} \label{BEYD}
\begin{itemize}
 \item
 A dessin $(X, {\cal D})$ defined on a topological (compact oriented) surface induces a canonical structure of Riemann surface on $X$ and a Belyi function
$$\beta : X \to \hat{\C}$$
such that ${\cal D}=\beta^{-1}([0,1])$.

\item
Conversely, if $$\beta :X \to \hat{\C}$$
is a  Belyi function defined on a (compact oriented) Riemann surface, ramified over at most the three points $0,1,\infty$, then ${\cal D}=\beta^{-1}([0,1])$ is a dessin on $X$.
\end{itemize}
\end{theorem}

\noindent In Theorem \ref{BEYD}, we suppose that the ramification values are $0,1, \infty$. It is not a restriction because, up to a Moebius transformation of  $\hat{\C}$, we can always  (without loss of generality) send  the white vertices of ${\cal D}$ onto $0$ and the black ones onto $1$. On the other hand, we remark that, although the structure of Riemann surface associated to the dessin by mean of a Belyi function  is unique, the Belyi function itself is not (it depends in particular on the position of the $*$-vertices of the triangulation and the homeomorphism from each butterfly onto $\S^2$).

\subsection{Coincidence of ${\cal C}_1$ and ${\cal C}_2$ for equilateral triangulations}\label{DEUXCOIN}
In this section, we will show that the Riemann structures ${\cal C}_1$  built in section \ref{SUBPO1}  and  ${\cal C}_2$ built in section \ref{SUBPO2}, induced by a dessin $(X, {\cal D})$ are identical if ${\cal D}$ is an equilateral triangulation.

\begin{theorem}\label{EQUST}
Let $X$ be a (compact oriented)  Riemann surface. The following assertions are equivalent~:
\begin{enumerate}
\item \label{EQUST1}
The conformal structure of $X$ is the ${\cal C}_2$-structure associated to the Belyi function obtained from a triangulation ${\cal T}$ on $X$.
\item \label{EQUST2}
The conformal structure of $X$ is the ${\cal C}_1$-structure associated to an equilateral Euclidean triangulation ${\cal T}$ on $X$.
\end{enumerate}
\end{theorem}

\noindent {\small {\bf Proof of Theorem \ref{EQUST}~-}
\begin{enumerate}
\item
Let us first suppose that the conformal structure of $X$ is the ${\cal C}_2$-structure associated to a Belyi function
$$\beta: X \to \hat{\C},$$
constructed from a triangulation ${\cal T}$ on $X$,
whose ramified values  belong to $\{0,1,\infty \}$. Let us consider the equilateral Euclidean triangulation ${\cal T}'$ defined on $\hat{\C}$ with vertices $\{0,1,\infty \}$. Then, ${\cal T}= \beta^{-1}{\cal T}'$ can be endowed with a metric such that each of its triangle $t$ is isometric to $\beta(t)$. This metric may have singularities at the vertices of ${\cal T}$. All triangles of ${\cal T}$ are isometric, and ${\cal T}$ is an equilateral Euclidean triangulation on $X$. Let $\{v_1,...,v_n\}$ denotes the set of vertices of ${\cal T}$. The covering $\beta$ induces a local isometry (and then a locally biholomorphic covering) from $X \backslash \{v_1,...,v_n\}$ onto $\hat{\C} \backslash \{0,1, \infty\}$.  Consequently, the conformal structure ${\cal C}_1$ of $X\backslash \{v_1,...,v_n\}$ (induced by the equilateral triangulation ${\cal T}$) coincides with the conformal structure ${\cal C}_2$ on $X\backslash \{v_1,...,v_n\}$ induces by $\beta$. Because $\{v_1,...,v_n\}$ ia a finite set, these structures coincide everywhere on $X$.

\item
Conversely, let us suppose that  the conformal structure of $X$ is the ${\cal C}_1$-structure associated to an equilateral Euclidean triangulation ${\cal T}$ on $X$. We will build a Riemann covering of $X$ over $\hat{\C}$ as follows~:

 \begin{figure}[ht!]
 \begin{center}
 \includegraphics[width=0.3\textwidth]{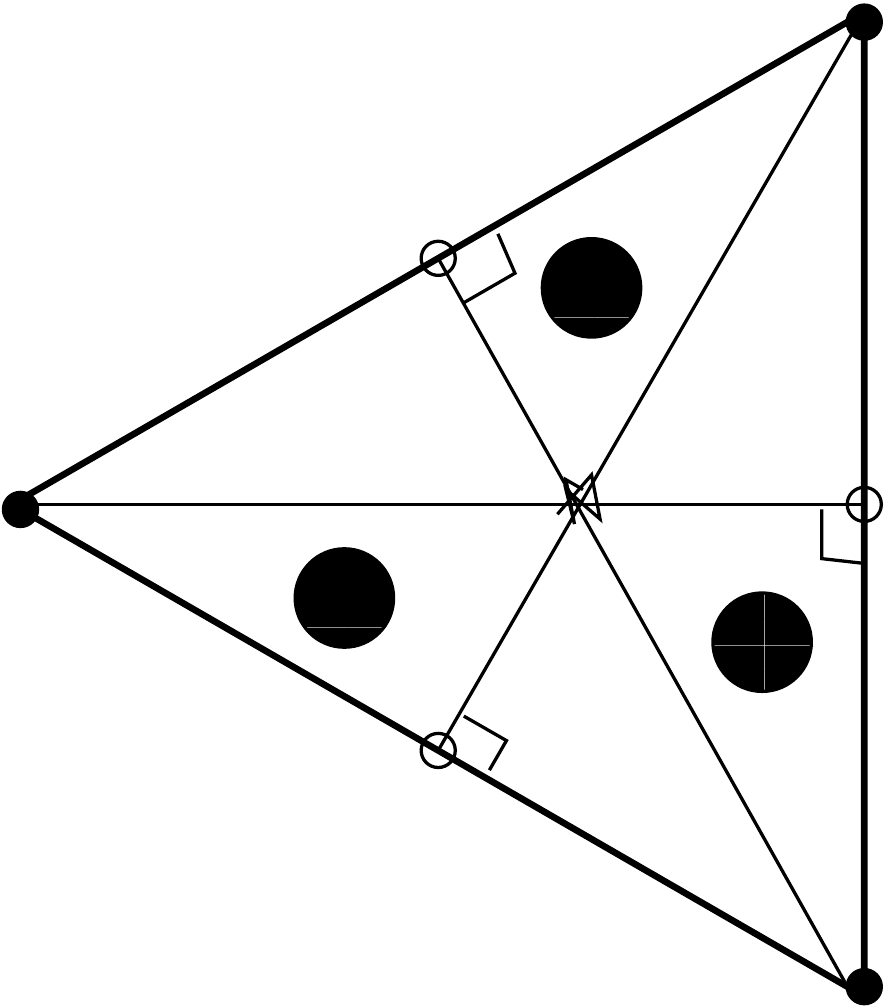}
 \end{center}
 \caption{A butterfly with right triangles}
 \end{figure}

  \begin{figure}[ht!]
 \begin{center}
 \includegraphics[width=0.6\textwidth]{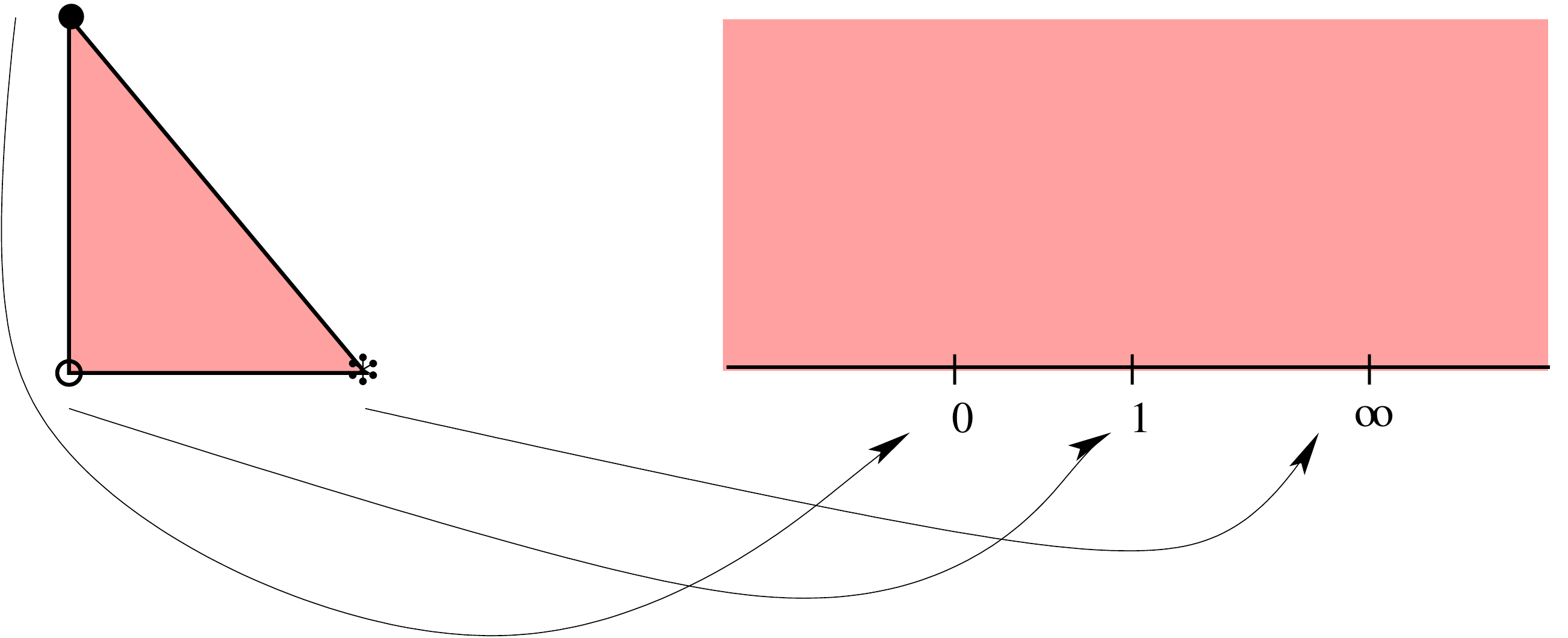}
 \end{center}
 \caption{Sending a right triangle onto $\C^+$}
 \end{figure}

  \begin{figure}[ht!]
 \begin{center}
 \includegraphics[width=0.6\textwidth]{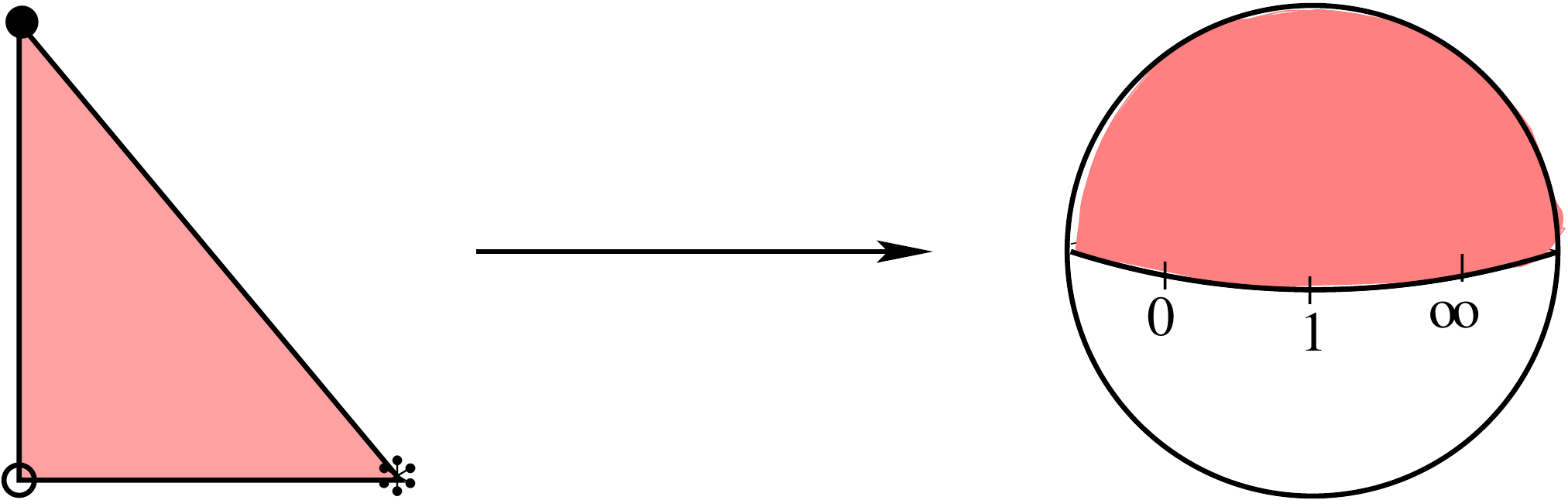}
 \end{center}
 \caption{Sending a right triangle onto the north hemisphere of $S^2 \simeq \hat{\C}$}
 \end{figure}

\begin{Myitemize}
\item First of all, we will build a "Euclidean butterfly decomposition" of each triangle $t$ of  ${\cal T}$ by building a new tricolored triangulation $\overline{{\cal T}}$ as follows~: We color the vertices of ${\cal T}$ in black.  By drawing the medians of each triangle $t$, we decompose each triangle of ${\cal T}$ in $6$ right triangles,  coloring in  white the vertices that are the  intersections of the edges of ${\cal T}$ with the medians and in $*$ the vertices that are the intersections of the medians. We get for each equilateral triangle $t \in {\cal T}$, $6$ triangles with angles of $\frac{\pi}{6}$, $\frac{\pi}{3}$, $\frac{\pi}{2}$.

\item
Then, we color alternatively each triangle of $\overline{{\cal T}}$ in black and white, denoting by $t^+$ the black ones and by $t^-$ the white ones. Each triangle becomes the union of $3$ butterflies $t^+ \cup t^-$. We get $6$ triangles with angles of $30$, $60$, $90$ degrees.

\item
 Now, we build a biholomorphism from each butterfly onto $\hat{\C}$ as follows.
 Using the Riemann mapping Theorem \ref{RIEMANNTH} and the (inverse of) the Riemann Christoffel transformation (Theorem \ref{RIEMANNCHRIS}), we build a holomorphic transformation $\beta^+$ of each triangle $t^+$ onto the upper plane $\hat{\C}^{+}$, sending the boundary of $t^+$ onto the boundary $\overline{\R}$, such that the black vertex is send onto $0 \in \R$, the white vertex onto $1 \in \R$ and the $*$-vertex onto $\infty$. By the same process, we build a holomorphic transformation $\beta^-$ of each triangle $t^-$ onto the lower plane $\hat{\C}^{-}$. Since $t^+$ and $t^-$ are isometric, $\beta^+$ and $\beta^-$ coincide on the common edge of $t^+$ and $t^-$. We obtain a continuous map $\beta^{\pm}$ from a butterfly onto $\hat{\C}$, that is biholomorphic except eventually on the common edge the $t^+$ and $t^-$. By the reminders of Section \ref{REMIND}, this transformation is a biholomorphism.

\item
We go on, by building a holomorphic covering of degree $3$  from each triangle $t \in {\cal T}$  over $\hat{\C}$
$$\beta_t : t \to \hat{\C},$$
of degree $3$
(ramified over $\infty$),
and then,  a holomorphic   covering
$$\beta : X \to \hat{\C},$$
(ramified over at most the three points $0,1, \infty$)
such that ${\cal T} = \beta^{-1}([0,1])$.
\end{Myitemize}

\noindent Finally, The covering $\beta$ is a Belyi  function, $(X, \beta)$ is a Belyi pair, from which we deduce that the conformal structure ${\cal C}_2$ associated to the dessin ${\cal T}$ {\it via} the Belyi function $\beta$ is nothing but ${\cal C}_1$.

\end{enumerate}
}

\section{Automorphisms of a dessin}
Classically, an {\it automorphism of graph} is a bijection of the set of vertices of the graph preserving the set of edges~: a pair of vertices is an edge if and only if its image by the bijection is also an edge. Let us now define an {\it automorphism of a dessin}. Although a purely combinatorial definition is possible, we prefer in our context a topological one.

\begin{definition} An {\rm automorphism of a dessin} $(X,{\cal D})$ is an isotopy class of homeomorphisms of $X$ preserving the graph ${\cal D}$ and the color of its vertices. We denote by $\mbox{\rm Aut}(X, {\cal D})$ the group of  automorphisms of the dessin $(X,{\cal D})$.
\end{definition}

\noindent Of course, two homeomorphisms of $X$ preserving ${\cal D}$ can induce the same automorphism of the graph ${\cal D}$.  However, one has the following important result \cite{jones2016dessins}, \cite{girondo2011introduction}~:

\begin{proposition} \label{HOMEO3}
In each isotopy class of automorphism of a dessin  $(X,{\cal D})$, there exists a unique (biholomorphic) automorphism of the Riemann surface $X$.
\end{proposition}

\noindent {\small {\bf Proof of Proposition \ref{HOMEO3}~-} Let ${\cal T}$ be the equilateral triangulation associated to ${\cal D}$. There is a unique isometry of $X$ endowed with the Euclidean structure associated to ${\cal T}$ on $X$, in an isotopy class of homeomorphism of the Riemann surface $X$ (endowed with the ${\cal C}_1={\cal C}_2$-conformal structure -~see Theorem \ref{EQUST}~-), preserving the graph ${\cal D}$ and the color of its vertices. This isometry induces a biholomorphic automorphism of $X$ preserving the ${\cal D}$ and the color of its vertices.
}

\vspace{10pt}

 \noindent Since ${\rm Aut}(X,{\cal D})$ is obviously finite, we deduce from Proposition \ref{HOMEO3}, Theorem \ref{UNIR} and Theorem \ref{SGSO}~:
 \begin{lemma} \label{LCONJ}
 Let $(X, {\cal D})$ be a dessin.  Then, there exists a canonical injective morphism of the group ${\rm Aut}(X,{\cal D})$ into the group ${\rm Aut}(X)$ of (biholomorphic) automorphisms of $X$, identifying canonically  ${\rm Aut}({\cal D})$  to a finite subgroup of ${\rm Aut}(X)$. In particular, ${\rm Aut}(X,{\cal D})$ is  conjugate to a finite subgroup of $SO(3)$.
 \end{lemma}

\section{Proof of Theorem \ref{THPR}}
Gathering together Theorem \ref{BEYD}, section \ref{KJOI9}, Lemma \ref{LCONJ},  Proposition \ref{SGSO21}, Theorem \ref{MOD}, Theorem \ref{THFI} and Theorem \ref{PKK}, we solve our initial problem~: Let $(X, \mathcal{D})$ be a dessin.  By section \ref{BF2}, we know that  $(X, {\cal D})$ admits a Belyi function and a canonical complex structure. We will distinguish the cases ${\bf g}_X \neq 1$ and ${\bf g}_X = 0$.

\subsection{The case ${\bf g}_X \neq 0$}
 \begin{Myitemize}
 \item
  If ${\bf g}_X > 1$,   $(X, {\cal D})$ admits a canonical Riemannian metric with constant Gaussian curvature $-1$, that is, the unique metric given by Corollary \ref{THFI}, invariant by the group of biholomorphisms of $X$.
  \item
  If ${\bf g}_X =1$,   $(X, {\cal D})$ admits, up to a scaling constant, a canonical flat Riemannian metric, that is, the  metric given by Theorem \ref{THFI}, invariant by the group of biholomorphisms of $X$.
  \end{Myitemize}

\subsection{The case ${\bf g}_X =0$} \label{SFSIT}
  If ${\bf g}_X = 0$, by Theorem \ref{MOD}, $X \simeq \hat{\C}$ admits a unique structure of Riemann surface, identified with $\hat{\C}$. We have seen in section \ref{KJOI9}  that there is no canonical Riemannian metric (invariant under Moebius transformations and of constant Gaussian curvature $1$) associated to this conformal structure. However the data of a dessin d'enfants ${\cal D}$ induces the data of the subgroup ${\rm Aut}(\hat{\C},{\cal D})$ of ${\rm Aut}(\hat{\C})$. That is why, we propose different approaches  introducing the automorphism group ${\rm Aut}(\hat{\C},{\cal D})$ to define  canonical Riemannian metrics on $\hat{\C}$. In each case, we use the fact that ${\rm Aut}(\hat{\C},{\cal D})$ acts on $X$ as a finite subgroup of ${\rm Aut}(\hat{\C})$.

   \begin{enumerate}
   \item
     Our first approach applies directly Proposition \ref{SGSO21}~:  $(\hat{\C}, {\cal D})$  can be canonically endowed with the Riemannian metric
   $$\tilde{g}= \frac{1}{{\rm card}({\rm Aut}(\hat{\C},{\cal D}))}\sum_{h \in {\rm Aut}(\hat{\C},{\cal D})}h^*(g_{\hat{\C}}),$$
      where $g_{\hat{\C}}$ is the standard metric of the round sphere of radius $1$. This metric $\tilde{g}$ obviously depends on ${\rm Aut}(\hat{\C},{\cal D})$,  and if ${\rm Aut}(\hat{\C},{\cal D})$ is included in $SO(3)$, $\tilde{g}$ coincides with the standard metric $g_{\hat{\C}}$ of the round sphere of radius $1$. We remark that the construction of this metric does not requires any restriction on the subgroup ${\rm Aut}(\hat{\C},{\cal D})$. However, the Gaussian curvature of this metric is not constant in general.

   \item
   Our second approach mimics the hyperbolic situation. We know that there exists a unique Riemannian metric of constant curvature $-1$ on $\H$, invariant by the group of biholomorphisms of ${\rm Aut}(\H)$ (see section \ref{KJOI77}). Although this strong property is no more true for $\hat{\C}$, we will build on $(\hat{\C}, {\cal D})$ a Riemannian metric that is invariant by the subgroup ${\rm Aut}(\hat{\C},{\cal D})$ of the group of Moebius transformations ${\rm Aut}(\hat{\C})$.
   \begin{itemize}
   \item
   Since ${\rm Aut}(\hat{\C},{\cal D})$ is a finite subgroup of $PSL(2,\C)$ we know by  the results of section \ref{FINSUB} that it is conjugate to a finite subgroup $K$ of $SO(3)$~: there exists $\varphi \in PSL(2,\C)$ such that
   $${\rm Aut}(\hat{\C},{\cal D}) = \varphi^{-1} K  \varphi.$$

   \item
   We define the Riemannian metric
   $$g=\phi^*g_{\hat{\C}},$$
   and we will prove that $g$ does not depend on $\phi$. Suppose that $\psi \in PSL(2,\C)$ satisfies
   $$\phi^{-1} K  \phi = \psi^{-1} K \psi.$$
   then,
   $$\psi \phi^{-1} K  \phi \psi^{-1},$$
   that is, $\phi \psi^{-1}$ belongs to the normalizer $N(K)$. We know by  Theorem \ref{NORMALIZER} that if $K$ is not cyclic, $N(K) \subset SO(3)$~: There exists $\alpha \in SO(3)$ such that $\psi=\alpha \varphi$.
   Then,
   $$\psi^*g_{\hat{\C}}= (\alpha \varphi)^*g_{\hat{\C}}= \varphi^*\alpha^*g_{\hat{\C}}=\varphi^*g_{\hat{\C}}.$$

   \item
   Now we prove that $g=\varphi^*g_{\C}$ is invariant by ${\rm Aut}(\hat{\C},{\cal D})$. As before, we know that there exists a finite subgroup $K$ of $SO(3)$ and $\varphi \in PSL(2,\C)$ such that
   $${\rm Aut}(\hat{\C},{\cal D})= \varphi^{-1} K \varphi.$$
   Let $\psi \in {\rm Aut}(\hat{\C},{\cal D})$. Let us compare $g$ and $\psi^*g$. There exists $\alpha \in K \subset SO(3)$ such that $\psi = \varphi^{-1} \alpha \varphi$. Then,
   $$\psi^*g=(\varphi^{-1} \alpha \varphi)^*g=\varphi^* \alpha^* (\varphi^{-1})^*\varphi^*g_{\C}=\varphi^* \alpha^*g_{\C}=\varphi^*g_{\C}=g.$$
   \end{itemize}
   \noindent Finally, $g$ is the solution of our problem.

   \end{enumerate}

\section{Addendum} \label{ADDENDUM}
We propose here two other constructions of a canonical Riemannian metric on a dessin $(X, {\cal D})$ when $X$ is a Riemann sphere. In these two last cases, the curvature of the metric is not constant in general.

\subsection{Figures of dessins on $\S^2$ invariant by a finite subgroup of $SO(3)$}
\begin{figure}[h]
\begin{center}
\includegraphics[height=5cm]{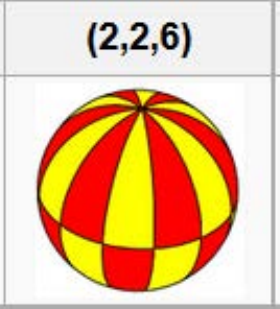}
\end{center}
\caption{Tesselations of the $2$-sphere, induced by the action of the diedral group ${\bf D}_6$, (From Wikipedia, Triangle group)
}
\label{diee}
\end{figure}

\begin{figure}
\begin{center}
    \includegraphics[height=4cm]{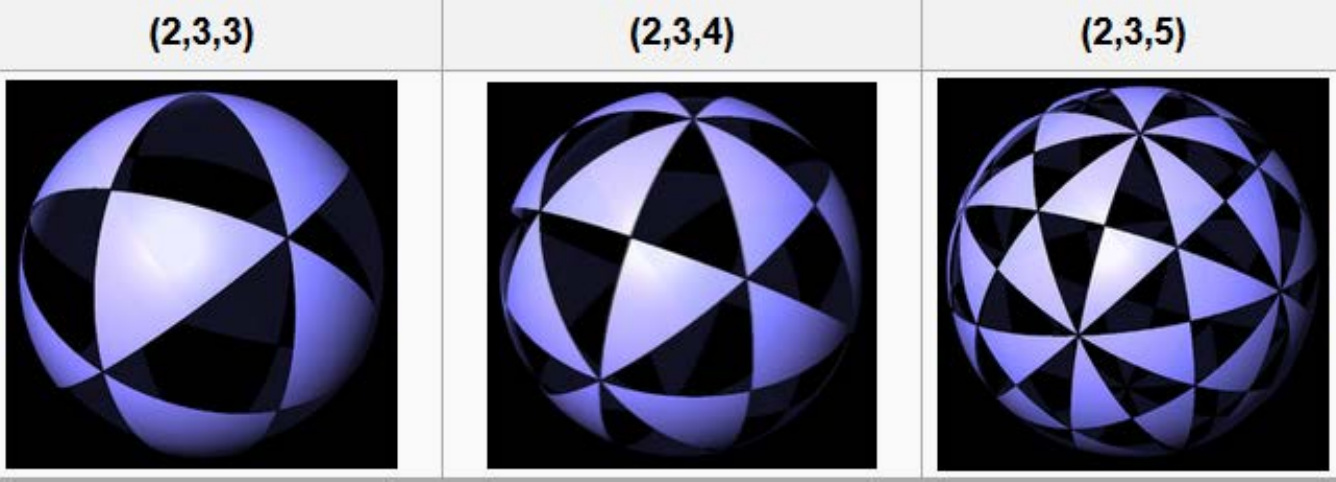}
\end{center}
\caption{Tesselations of the $2$-sphere, induced by the actions of ${\cal A}_4$, ${\cal S}_4$, ${\cal A}_5$, (From Wikipedia, Triangle group)
 }
\label{diee2}
\end{figure}

Figures  \ref{diee} and \ref{diee2} shows  the tessellations on $\S^2$ invariant by $\bf{D}_n$ (from Wikipedia, Triangle group) and $\bf{A}_4$, $\bf{S}_4$, $\bf{A}_5$ (from Wikipedia, Triangle group -~by Jeff Weeks~-).

\subsection{A third construction}
    A third method uses the following proposition~:
    \begin{proposition} \label{SGSO22}
  Let $K$ be a finite subgroup of ${\rm Aut}(\hat{\C})$. Then, $\C^2$ admits a canonical Hermitian metric $(.,.)_K$ invariant by $K$, whose real part $\tilde{\tilde{g}}_K$ induces on the Riemann sphere $\hat{\C}$ a Riemannian metric that coincide with the metric $g_{\hat{\C}}$  of the round sphere if $K$ is a subgroup of $SO(3)$.
 \end{proposition}

\noindent Remark that
Proposition \ref{SGSO22} shows that $K$ is conjugate to a finite subgroup of $SU(2)/\pm 1 \simeq SO(3)$.

 \vspace{10pt}

\noindent {\footnotesize {\bf Proof of Proposition \ref{SGSO22}~-}
 We consider the sequence of canonical embeddings
$$\S^2 \subset \C \times \R \subset \C \times \R^{\C} \simeq \C \times \C=\C^2,$$
where $\S^2  \subset \C \times \R$ is the standard totally umbilic isometric embedding of the round sphere of radius $1$ in $\R^3 \simeq \R^2 \times \R \simeq \C \times \R$. We deduce the standard isometric embedding
\begin{equation} \label{KAAR}
 \S^2  \to \S^3 \subset \C^2,
\end{equation}
where $\S^3$ is the standard round sphere of radius $1$ in $\R^4 \simeq \C^2$.
 Let us take the standard Hermitian scalar product $(.|.)$ on $\C^2$. If $K$ is any subgroup of ${\rm Aut}(\hat{\C}) \simeq PSL(2,\C)$, one can define a new Hermitian scalar product $(.|.)_K$ on $\C^2$ as follows~: For all $u$, $v$ in $\C^2$,
$$(u|v)_K= \frac{1}{{\rm card}(K)}\sum_{h \in K}(h(u)|h(v)).$$
(This new Hermitian scalar product is obviously $H$ invariant.)
Classically, let us decompose $(.|.)_K$ in its real part and its imaginary part~:
$$(.|.)_K= \tilde{\tilde{g}}_K  + \imath \Omega_K,$$
where $\tilde{\tilde{g}}_K$ is a Riemannian metric and $\Omega_K$ a symplectic form.
 Using \ref{KAAR}, we build  on $\S^2$ the Riemannian metric $\overline{g}_K=i^{*}\tilde{\tilde{g}}_K$. {\it A priori}, $\overline{g}_K$ is {\it not} a Riemannian metric with constant curvature $1$.
}
 \vspace{10pt}

  \subsection{A fourth construction}
   Our fourth method uses the study of the orbits of ${\rm Aut}(X,{\cal D})$, using the classification of the  subgroups of $PSL(2,\C)$ and its consequence~: If ${\rm Aut}(X,{\cal D})$ is not isomorphic to $\Z_n$, it contains exactly $3$ orbits  $O_1, O_2, O_3$.

    \begin{itemize}
      \item
      If ${\rm card (O_1)} <  {\rm card (O_2)} < {\rm card (O_3)}$, we associate
      to each triplet $(a_1, a_2, a_3) \in O_1 \times O_2 \times O_3$,  the unique Moebius transformation $h_{a_1a_2a_3}$ that sends $0$ to $a_1$, $1$ to $a_2$ and $2$ to $a_3$. Then, we build the metric
   $$g= \frac{1}{{\rm card}\{a_1a_2a_3\}}\sum_{a_1a_2a_3}h^{*}_{a_1a_2a_3}g_{\hat{\C}},$$
   where the sum is over all possible triplets $a_1a_2a_3$ built as before.

    \item
    If  ${\rm card (O_1)} =  {\rm card (O_2)} < {\rm card (O_3)}$ (it is the dihedral situation), we cannot distinguish the orbits with the same cardinality, so with the same notations,  we associate
      to each triplet $(a_1, a_2, a_3) \in O_1 \times O_2 \times O_3$,  the unique Moebius transformation $h_{a_1a_2a_3}$ that sends $0$ to $a_1$, $1$ to $a_2$ and $2$ to $a_3$, and the unique Moebius transformation $k_{a_1a_2a_3}$ that sends $0$ to $a_2$, $1$ to $a_1$ and $2$ to $a_3$, . Then we build the metric
   $$g= \frac{1}{2{\rm card}\{a_1a_2a_3\}}\sum_{a_1a_2a_3}h^{*}_{a_1a_2a_3}g_{\hat{\C}}+
   \frac{1}{2{\rm card}\{a_1a_2a_3\}}\sum_{a_1a_2a_3}k^{*}_{a_1a_2a_3}g_{\hat{\C}}.$$

    \item
    If  ${\rm card (O_1)} <  {\rm card (O_2)} = {\rm card (O_3)}$ (it is the tetrahedral situation), an analogous process builds  a canonical metric.
    \end{itemize}

    \noindent Finally, if ${\rm Aut}(X, {\cal D})$ is isomorphic to $\Z_n$, we can arbitrarily endow $\hat{\C}$ with the standard metric $g_{\C}$.

  \vspace{10pt}

  \noindent The reader can produce other Riemannian metrics of this type on $\hat{\C}$, playing for instance with the fixed points of the elements of ${\rm Aut}(X, {\cal D})$.

\newpage

%
%

\vspace{3cm}

\newpage

\bibliographystyle{plain}


\begin{thebibliography}{1}
\bibitem{bost1992introduction} J.-B. Bost. {\it Introduction to compact Riemann surfaces, Jacobians, and Abelian varieties}.  From number theory to physics, pp. $64-211$. Springer, $1992$.

\bibitem{cheltsov2016} I. Cheltsov and C. Shramov. {\it Cremona groups and the icosahedron}.
  Monographs and Research Notes in Mathematics,
  $2016$, CRC Press.

\bibitem{dai2008variational} J. Dai and X. D. Gu and F. Luo.
  {\it Variational principles for discrete surfaces}.
  Volume $4$,
  $2008$,
  International Press of Boston Incorporated.

  \bibitem{donaldson2011riemann} S. Donaldson. {\it Riemann surfaces}.
 $2011$. Oxford University Press.

 \bibitem{earp1997} R. Earp and E. Toubiana.
  {\it Introduction  a la geom\'{e}trie hyperbolique et aux surfaces de Riemann}.
  $1997$.
  Diderot \'{e}diteur arts et sciences.

  \bibitem{girondo2011introduction} E. Girondo and G. Gonz{\'a}lez-Diez.
  {\it Introduction to compact Riemann surfaces and dessins d’enfants},
  Volume $79$,
  $2011$.
  Cambridge University Press.

\bibitem{grothendieck1984esquisse} A. Grothendieck.
  {\it Esquisse d'un programme}.
  Preprint, Montpellier,
  $1984$.

\bibitem{gu2008computational} X. D. Gu and S.T. Yau.
  {\it Computational conformal geometry}. $2008$,
  International Press Somerville, Mass, USA.


\bibitem{guillarmoumoroianu} C. Guillarmou and S. Moroianu.
  {\it Unpublished, web site of C. Guillarmou}.

\bibitem{jones2016dessins} G. Jones, A. Gareth  and J. Wolfart.
  {\it Dessins d'enfants on Riemann surfaces}.
  $2016$, Springer.

\bibitem{lando2013graphs} S.K. Lando and A. K. Zvonkin.
  {\it Graphs on surfaces and their applications}.
  Volume $141$,
  $2013$,
  Springer Science, Business Media

\bibitem{malgoire12cartes} J. Malgoire and C. Voisin.
  {\it Cartes cellulaires}.
  Cahiers Math. Montpellier,
  volume $12$.

\bibitem{REYSSAT} E. Reyssat {\it Quelques aspects des surfaces de Riemann}.
 $1989$, Progress in Math. $77$, Birkhauser.


\bibitem{shabat1990drawing} G. B. Shabat and V.A. Voevodsky.
  title={Drawing curves over number fields}.
   The Grothendieck Festschrift,
  $199-227$,
  $1990$,
  Springer.


\bibitem{voisin1980cartes} Ch. Voisin and J. Malgoire.
  {\it Cartes topologiques infinies et rev{\^e}tements ramifi{\'e}s de la sph{\`e}re}.
 $1980$,
  Univ. des Sciences et Techniques du Languedoc, UER de Math{\'e}matiques.


\end{thebibliography}

\Addresses

\end{document}